\documentclass[a4paper,11pt]{amsart}
\usepackage{amssymb}
\usepackage{amscd}

\def\Box{\square}

\def\frak{\mathfrak}
\def\Bbb{\mathbb}
\def\Cal{\mathcal}

\def\bbW{{\mathbb{W}}}
\def\bbX{{\mathbb{X}}}
\def\bbY{{\mathbb{Y}}}
\def\bbZ{{\mathbb{Z}}}

\newcommand{\hook}{\raisebox{-0.35ex}{\makebox[0.6em][r]
{\scriptsize $-$}}\hspace{-0.15em}\raisebox{0.25ex}{\makebox[0.4em][l]{\tiny
 $|$}}}

\let\hash=\sharp

\let\d=\delta
\let\e=\varepsilon
\let\f=\varphi
\let\ii=\iota

\let\w=\omega
\let\r=\rho

\def\uw{\underline{\wedge}}

\def\fsl{{\frak{sl}}}

\def\LOT{{\rm LOT}}

\newcommand{\miniX}{\mbox{\boldmath{$\scriptstyle{X}$}}}

\newcommand{\LX}{{\Cal L}_{\!\miniX}}
\newcommand{\KX}{{\Cal K}_{\!\miniX}}

\newcommand{\NX}{{\mbox{\boldmath$ \nabla$}}_{\!\miniX}}

\let\Euler=\NX
\newcommand{\IT}[1]{{\rm(}{\it{\!#1}}{\rm)}}

\newcommand{\newc}{\newcommand}

\newcommand{\Ker}{\operatorname{Ker}}

\let\ccdot\cdot
\def\cdot{\hbox to 2.5pt{\hss$\ccdot$\hss}}

\newcommand{\om}{\omega}
\renewcommand{\phi}{\varphi}
\newcommand{\ph}{\varphi}

\newcommand{\si}{\sigma}

\newc{\aI}{\mbox{\boldmath{$ I$}}}
\newc{\aR}{\mbox{\boldmath{$ R$}}}
\newc{\aDeR}{\mbox{\boldmath{$ U$}}_B{}^P{}_C{}^Q}
\newc{\al}{\mbox{\boldmath$ \Delta$}}             
\newc{\nda}{\mbox{\boldmath$ \nabla$}}
\newc{\ad}{\mbox{\boldmath$ d$}}
\newc{\da}{\mbox{\boldmath$ \delta$}}
\newc{\aK}{\mbox{\boldmath{$ K$}}}
\newc{\aL}{\mbox{\boldmath{$ L$}}}

\newtheorem{theorem}{Theorem}[section]
\newtheorem{lemma}[theorem]{Lemma}
\newtheorem{proposition}[theorem]{Proposition}

\newcommand{\cC}{{\Cal C}}

\newcommand{\cH}{{\Cal H}}
\newcommand{\cg}{{\Cal G}}
\newcommand{\cN}{{\Cal N}}
\newcommand{\cV}{{\Cal V}}

\newcommand{\ce}{{\Cal E}}

\newcommand{\cq}{{\Cal Q}}

\newcommand{\cL}{{\Cal L}}
\newcommand{\cT}{{\Cal T}}
\newcommand{\cR}{{\Cal R}}

\newcommand{\ct}{{\Cal T}}

\newcommand{\bV}{{\Bbb V}}
\newcommand{\bU}{{\Bbb U}}
\newcommand{\bL}{{\Bbb L}}
\newcommand{\bG}{{\Bbb G}}

\newcommand{\nd}{\nabla}

\newcommand{\Rho}{{\mbox{\sf P}}}
\newcommand{\Up}{\Upsilon}

\newcommand{\End}{\operatorname{End}}

\newcommand{\Ric}{\operatorname{Ric}}

\newcommand{\horisys}[2]{(#1,#2)}
\newcommand{\vertsys}[2]{\left(\begin{array}{c}#1\\#2\end{array}\right)}

\newcommand{\Q}{\mbox{$
\begin{picture}(9,8)(1.6,0.15)
\put(1,0.2){\mbox{$ \mbox{\boldmath$ Q$} \hspace{-7.8pt} /$}}
\end{picture}$}}

\newcommand{\fb}  {\mbox{$                                      
\begin{picture}(9,8)(1.6,0.15)
\put(1,0.2){\mbox{$ \Box \hspace{-7.8pt} /$}}
\end{picture}$}} 

\newcommand{\afl}{\mbox{$
\begin{picture}(9,8)(1.6,0.15)
\put(1,0.2){\mbox{$ \al \hspace{-7.8pt} /$}}
\end{picture}$}}

\newcommand{\fD}{\mbox{$
\begin{picture}(9,8)(1.6,0.15)
\put(1,0.2){\mbox{$ D \hspace{-7.8pt} /$}}
\end{picture}$}}

\newcommand{\afD}                         
{\mbox{$
\begin{picture}(9,8)(1.6,0.15)
\put(1,0.2){\mbox{$ \D \hspace{-7.8pt} /$}}
\end{picture}$}}

\newcommand{\fl}{\mbox{$                                     
\begin{picture}(9,8)(1.6,0.15)
\put(1,0.2){\mbox{$ \Delta \hspace{-7.8pt} /$}}
\end{picture}$}}

\newcommand{\act}{\mbox{\boldmath{$\ct$}}}

\newcommand{\td}{\tilde{d}}
\newcommand{\dt}{\tilde{\d}}

\newcommand{\bK}{{\Bbb K}}

\newcommand{\nn}[1]{(\ref{#1})}

\newcommand{\D}{\mbox{\boldmath{$ D$}}}

\newcommand{\bF}{\mbox{{$\Bbb F$}}}
\newcommand{\bT}{\mbox{{$\Bbb T$}}}

\newcommand{\X}{\mbox{\boldmath{$ X$}}}

\newcommand{\sX}{\mbox{\scriptsize\boldmath{$X$}}}        
\newcommand{\h}{\mbox{\boldmath{$ h$}}}
\newcommand{\bg}{\mbox{\boldmath{$ g$}}}

\newcommand{\cce}{\tilde{\ce}}                          
\newcommand{\aM}{\tilde{M}}

\newcommand{\tU}{\tilde{U}}

\newcommand{\sbg}{\mbox{\scriptsize\boldmath{$g$}}}

\let\m=\mu

\newcommand{\V}{{\mbox{\sf P}}}                   
\newcommand{\J}{{\mbox{\sf J}}}

\newc{\strutdd}{\rule{0mm}{5mm}}

\newcommand{\rpl}                         
{\mbox{$
\begin{picture}(12.7,8)(-.5,-1)
\put(0,0.2){$+$}
\put(4.2,2.8){\oval(8,8)[r]}
\end{picture}$}}

\newcommand{\lpl}                         
{\mbox{$
\begin{picture}(12.7,8)(-.5,-1)
\put(2,0.2){$+$}
\put(6.2,2.8){\oval(8,8)[l]}
\end{picture}$}}

\newcommand{\upl}                         
{\mbox{$                                  
\begin{picture}(12.7,8)(-.5,-1)
\put(2,0.2){$+$}
\put(6.2,2.8){\oval(8,8)[t]}
\end{picture}$}}

\newcommand{\bGuk}                         
{\mbox{$
\begin{picture}(15,25)(-.5,15)
\put(2,28){$\ce^{k-1}$}
\put(0,16){\mbox{$\upl$}}
\put(2,4){$\ce^k$}        
\end{picture}$}}

\newcommand{\bGdk}                         
{\mbox{$
\begin{picture}(15,25)(-.5,15)
\put(2,28){$\ce_{k}$}
\put(0,15){\mbox{$\upl$}}
\put(2,4){$\ce_{k-1}$}        
\end{picture}$}}

\newcommand{\bGdkone}                         
{\mbox{$
\begin{picture}(20,25)(-.5,15)
\put(2,28){$\ce_{k+1}$}
\put(0,15){\mbox{$\upl$}}
\put(2,4){$\ce_{k}$}        
\end{picture}$}}

\newcommand{\cbGuk}{\mbox{$\begin{picture}(50,35)(-.5,15)
                            \put(0,16){$\cg^k=\bGuk$}\end{picture}$}}

\newcommand{\cbGdk}{\mbox{$\begin{picture}(50,35)(-.5,15)
                            \put(0,16){$\cg_k=\bGdk$}\end{picture}$}}

\newcommand{\cbGdkone}{\mbox{$\begin{picture}(60,28)(-.5,15)
                            \put(0,16){$\cg_{k+1}=\bGdkone$}\end{picture}$}} 

\usepackage{ifthen}

\newc{\tensor}[1]{#1}
\newc{\Mvariable}[1]{\mbox{#1}}
\newc{\down}[1]{{}_{
\ifthenelse{\equal{#1}{;}}{|}{#1}}}
\newc{\up}[1]{{}^{#1}}
\newc{\C}{C}

\newc{\JulyStrut}{\rule{0mm}{6mm}}
\newc{\midtenPan}{\mbox{\sf S}}
\newc{\midten}{\mbox{\sf T}}
\newc{\midtenEi}{\mbox{\sf U}}
\newc{\ATen}{\mbox{\sf E}}
\newc{\BTen}{\mbox{\sf F}}
\newc{\CTen}{\mbox{\sf G}}

\def\sideremark#1{\ifvmode\leavevmode\fi\vadjust{\vbox to0pt{\vss
 \hbox to 0pt{\hskip\hsize\hskip1em
 \vbox{\hsize3cm\tiny\raggedright\pretolerance10000
 \noindent #1\hfill}\hss}\vbox to8pt{\vfil}\vss}}}%

\begin{document}
\renewcommand{\today}{}

\title{Conformal de Rham Hodge theory and operators generalising the
Q-curvature} \author{A. Rod Gover}

\address{Department of Mathematics\\
  The University of Auckland\\
  Private Bag 92019\\
  Auckland 1\\
  New Zealand} \email{gover@math.auckland.ac.nz}

\vspace{10pt}

\renewcommand{\arraystretch}{1}
\maketitle
\renewcommand{\arraystretch}{1.5}

\pagestyle{myheadings}
\markboth{Gover}{Conformal Hodge theory and a generalisation of
Q-curvature}

\begin{abstract}
  We look at several problems in even dimensional conformal geometry
  based around the de Rham complex. A leading and motivating problem
  is to find a conformally invariant replacement for the usual de Rham
  harmonics.  An obviously related problem is to find, for each order
  of differential form bundle, a ``gauge'' operator which completes
  the exterior derivative to a system which is both elliptically
  coercive and conformally invariant. Treating these issues involves
  constructing a family of new operators which, on the one hand,
  generalise Branson's celebrated Q-curvature and, on the other hand,
  compose with the exterior derivative and its formal adjoint to give
  operators on differential forms which generalise the critical
  conformal power of the Laplacian of Graham-Jenne-Mason-Sparling.
  We prove here that, like the critical conformal Laplacians, these 
  conformally invariant operators are not strongly invariant. 
 The construction draws heavily on the ambient metric of
  Fefferman-Graham and its relationship to the conformal tractor
  connection and exploring this relationship will be a central theme
  of the lectures.
\end{abstract}
\thanks{ARG gratefully acknowledges support from the Royal Society of
New Zealand via a Marsden Grant (grant no.\ 02-UOA-108). The author is
a New Zealand Institute of Mathematics Maclaurin Fellow.}

These notes draw on recent collaborative work with Tom Branson. There
is also significant input from recent joint work with Andi \v Cap and
Larry Peterson. Among the results that are completely new here is
Proposition \ref{notstrong} which proves that the conformally
invariant differential operators between forms, that we construct
here, have the curious property that thay are not strongly invariant.
Also it is shown, in the final section, that the ``$Q$-operators'',
which were first developed in \cite{deRham}, can be recovered by a
polynomial continuation argument that parallels and generalises
Branson's original construction of the $Q$-curvature.  These notes were
presented as a series of three lectures at the 24${}^{\rm th}$ Winter
School on Geometry and Physics, Srn\'{\i} Czech Republic, January
2004.


\section{Lecture 1 -- some problems related to the de Rham complex}

The de Rham complex and its cohomology are among the most
fundamental of tools for relating local differential geometric
information to global topology. Over these lectures we shall explore
the de Rham complex and related issues in the setting of conformal geometry.

On a smooth $n$-manifold $M$, let us write $\ce$ or $\ce^0$ for
$C^\infty(M)$ and $\ce^k$ for the space of k-forms, i.e., the smooth
sections of the k${}^{\rm th}$ exterior power of the cotangent bundle
$\Lambda^kT^*M$.  Recall that the exterior derivative on functions
takes values in $T^*M$ and is defined by $df(v)=vf$ where we view the
smooth tangent vector field $v$ as a derivation (so in terms of local
coordinates $x^i$ we have $vf=\sum v^i\partial f/\partial x^i$). This
is extended to a differential operator $d:\ce^k\to \ce^{k+1}$ by
requiring $d^2 f:=d(d f)=0$, for $f\in \ce$, and the Leibniz rule $d f
w = df\wedge w + fd w$, $w\in\ce^k$. It follows that $d^2$ vanishes on k-forms and
so we obtain the de Rham complex,
$$
\ce^0\stackrel{d}{\to} \ce^1\stackrel{d}{\to} \cdots \stackrel{d}{\to}\ce^{n}.
$$

We are interested in the additional operators between the form bundles
that arise when $M$ is equipped with a conformal structure. Recall
that a {\em conformal structure} is an equivalence class of metrics
$[g]$ where two metrics are equivalent if they are related by
multiplication by a smooth positive function, i.e., $\widehat{g}\sim
g$ means there is $\w\in \ce$ such that $\widehat{g}=e^{2\w}g$.
We may equivalently view the conformal class as being
given by a smooth ray subbundle $\cq\subset S^2T^*M$, whose fibre at $x$ is
formed by the values of $g_x$ for all metrics $g$ in the conformal
class. By construction, $\cq$ has fibre $\Bbb R_+$ and the metrics in
the conformal class are in bijective correspondence with smooth
sections of $\cq$.

The bundle $\pi:\cq\to M$ is a principal bundle with structure group
$\Bbb R_+$, and we denote by $E[w]$ the line bundle induced from the
representation of $\Bbb R_+$ on $\Bbb R$ given by $s\mapsto s^{-w/2}$.
Sections of $E[w]$ are called a {\em conformal densities of weight
$w$} and may be identified with functions on $\cq$ that are
homogeneous of degree $w$, i.e., $f(s^2 g_x,x)=s^w f(g_x,x)$ for any
$s\in \Bbb R_+$.  We write $\ce[w]$ for the space of sections
of the bundle and, for example, $\ce^k[w]$ is the space of sections of
$(\Lambda^k T^*M)\otimes E[w]$. (Here and elsewhere all sections are
taken to be smooth.)

There is a tautological function $\bg$ on $\cq$ taking values in
$\ce_{(ab)}$.  It is the function which assigns to the point
$(g_x,x)\in \cq$ the metric $g_x$ at $x$.  This is homogeneous of
degree 2 since $\bg (s^2 g_x,x) =s^2 g_x$. If $\si$ is any positive
function on $\cq$ homogeneous of degree $1$ then $\si^{-2} \bg$ is
independent of the action of $\Bbb R_+$ on the fibres of $\cq$, and so
$\si^{-2} \bg$ descends to give a metric from the conformal
class. Thus $\bg$ determines and is equivalent to a canonical section
of $\ce_{ab}[2]$ (called the conformal metric) that we also denote
$\bg$ (or $\bg_{ab}$). Then, for $\si\in \ce_+[1]$, $\si^{-2}\bg$ is a
metric from the conformal class and we term $\si$ a {\em conformal scale}.
We will use the conformal metric to raise and lower indices.

Recall that the Levi Civita connection is the unique torsion free
connection on tensor bundles which preserves a given metric. So on a
conformal manifold a choice of conformal scale $\si$ determines a Levi
Civita connection that we will denote $\nabla$. The scale $\si$ also
determines a connection (that we also denote $\nabla$ and term the
Levi Civita connection) on densities by the formula $\nabla \mu= \si^w
d \si^{-w}\mu$, for $\mu\in \ce[w]$.  (Note that $\si^w d \si^{-w}\mu$
means $\si^w (d (\si^{-w}\mu))$. The default is that in such
expressions all symbols, except the one at the extreme right,
are to be interpreted as operators and parentheses are usually omitted.) 
 For $g=\si^{-2}\bg$ the {\em conformal
rescaling} $g\mapsto \widehat{g}=e^{2\w}g$ corresponds to $\si\mapsto
\hat{\si}=e^{-\w}\si$ and so it follows at once that
\begin{equation}\label{densconntrans}
\widehat{\nd} \mu=\nd \mu + w\Up \mu, 
\end{equation}
 where $\widehat{\nd}$ is the connection for $\hat{\si}$ and
$\Up:=d\w$. It is similarly easy to show that, for example on 1-forms,
the Levi Civita connection transforms conformally according to
\begin{equation}\label{conntrans}
\widehat{\nabla}_a u_b=\nd_a u_b -\Up_a u_b -\Up_b u_a + \bg_{ab} \Up^c u_c,
\end{equation}
 where abstract indices are used in an obvious way and the inverse
of $\bg$ is used to raise the index on $\Up^c$.


To simplify our subsequent discussion let us assume that $M$ is
connected, compact and orientable, and that the conformal structure is
Riemannian (i.e., $g\in [g]$ has Riemannian signature). Via the
conformal metric, the bundle of volume densities can be canonically
identified with $ E[-n]$ and so the Hodge star operator (for each
metric from the conformal class) induces a conformally invariant
isomorphism that we shall also term the Hodge star operator: $\star:
\ce^k \cong \ce^{n-k}[n-2k]$. Let us write $\ce_{n-k}$ as an
alternative notation for the image space here, so we have $\star:
\ce^k \cong \ce_{n-k}$ and more generally $\ce_k[w]:=\ce^k[w+2k-n]$.
This notation is suggested by the duality between the section spaces
$\ce^k$ and $\ce_k $. For $\phi\in \ce^k$ and $\psi\in \ce_k $, there
is the natural conformally invariant global pairing
$$
\phi,\psi \mapsto \langle \phi,\psi \rangle :=\int_M \phi\cdot \psi\,
d\mu_{\sbg}= \int_M \phi \wedge \star \psi ,
$$
where $\ph\cdot \psi\in\ce[-n]$
denotes a complete contraction between $\phi $ and $\psi$. 

In even dimensions $\ce^{n/2}=\ce_{n/2}$ and the de Rham
complex may be written in the more symmetric form
$$ 
\ce^0\stackrel{d}{\to}\cdots \stackrel{d}{\to}
\ce^{n/2-2}\stackrel{d}{\to}\ce^{n/2-1}\stackrel{d}{\to}
\ce^{n/2}\stackrel{\star
d}{\to}\ce_{n/2-1}\stackrel{\delta}{\to}\ce_{n/2-2}
\stackrel{\delta}{\to} \cdots\stackrel{\delta}{\to}\ce_0,
$$ where $\d$ is the composition $\star d \star$. (Or one could
 alternatively replace the $\stackrel{d}{\to}
\ce^{n/2}\stackrel{\star
d}{\to}$ with $ \stackrel{\star d}{\to}
 \ce^{n/2}\stackrel{\d}{\to}$.) Of course these operators are all
 conformally invariant since the exterior derivative is well defined
 on any smooth manifold and the conformal structure is used here only
 to give the isomorphisms $\star:\ce^k\to \ce_{n-k}$ and
 $\star:\ce_k\to \ce^{n-k}$.  Since $\star$ maps $\ce^{n/2}$ to itself
 we also have the conformally invariant operator $\d: \ce^{n/2}\to
 \ce_{n/2-1}$. This does not annihilate the exact forms; in fact in
 terms of the global pairing introduced above we have $\langle \f, \d
 d\f \rangle = \langle d\f, d\f \rangle $ and so, in the compact
 setting, $\d d\f =0$ implies $d\f=0$.  Rather the composition $\d
 d:\ce^{n/2-1}\to \ce_{n/2-1}$ is the well known conformally invariant
 {\em Maxwell operator}. (In dimension 4 and Lorentzian signature this
 gives the equations of electromagnetism.)  Thus we have the Maxwell
 {\em detour complex},
 $$ 
\ce^0\stackrel{d}{\to}\cdots \stackrel{d}{\to}
\ce^{n/2-2}\stackrel{d}{\to}\ce^{n/2-1}\stackrel{\d d}{\to}
\ce_{n/2-1}\stackrel{\delta}{\to}\ce_{n/2-2}
\stackrel{\delta}{\to} \cdots\stackrel{\delta}{\to}\ce_0.
$$ This really is symmetric since, for each $k$, the operator
$\d:\ce_{k+1}\to \ce_k$ is (up to a sign) the formal adjoint of
$d:\ce^k\to \ce^{k+1}$. It is convenient to absorb this sign and
redefine $\d$ to be exactly the formal adjoint, i.e., so that $\langle
\phi, \d \psi \rangle= \langle d\phi, \psi \rangle$ for $\phi\in \ce^k
$ and $\psi\in \ce_{k+1}$. The Maxwell operator and this detour
complex are a feature of even dimensional conformal geometry that we
wish to study and generalise. There are not analogues in odd
dimensions and so, for the remainder of this lecture (except where
otherwise indicated), let us suppose $n\geq 4$ is even.

As a point on notation. We will use $\ii(\cdot)$ and $\e(\cdot)$ as the
notation for interior and exterior multiplication by 1-forms on
differential forms. For a 1-form $u$ and a $k$-form $v$ the
conventions are
$$
(\e(u)v)_{a_0\cdots a_k}=(k+1)u_{a_0}v_{a_1\cdots a_k},\quad \mbox{and} 
\quad (\ii(u)v)_{a_2\cdots a_k}=u^{a_1}v_{a_1a_2\cdots a_k}.
$$ 
Here, and below, sequentially labelled indices are implicitly skewed over.

\subsection{The problems} We are now set to state and consider a series of 
fundamental problems 
concerning the de Rham complex on conformal manifolds.\\

\noindent{\bf Problem 1: Conformal Hodge theory.} Let us write
$H^k(M)$ for the $k^{\rm th}$ cohomology space of the de Rham complex.
If $M$ is equipped with a Riemannian metric $g$, then de Rham Hodge
theory exhibits an isomorphism between $H^k(M)$ and the space of {\em harmonics} $
\cH^k(M)$. The latter is the null space of the form Laplacian $\fl = \d d +d\d $ on $
k$-forms or, alternatively, it is recovered by
$$
\cH^k(M)=\cN(d:\ce^k\to\ce^{k+1})\cap \cN(\d : \ce^k \to
\ce^{k-1}).
$$

Here $\d$ is not, in general, the conformally invariant operator
described above but rather just the Riemannian formal adjoint of $d$.
For $u=u_{a_1\cdots a_k}$ a $k$-form we have 
$$
(du)_{a_0\cdots a_k}=(k+1)\nabla_{a_0}u_{a_1\cdots a_k},
$$ 
since the Levi Civita connection is torsion free. Thus $\d$ is given by 
$$
v_{a_0\cdots a_k}\mapsto  - \nd^{a_0} v_{a_0\cdots a_k}.
$$
This agrees with the conformally invariant $\d$ only if $v\in
\ce_k$; since only then is the integration by parts a conformally
invariant operation. (A point on notation: we will write $d u$ to mean
$(k+1)\nabla_{a_0}u_{a_1\cdots a_k}$ and $\d v$ to mean $-\nd^{a_0}
v_{a_0\cdots a_k}$ even when the density weights of $u$ and $v$ are
such that these are not conformally invariant.) Note that
$\ce^k=\ce_k\otimes\ce[n-2k]$. Thus from the Leibniz rule for $\nabla$
and \nn{densconntrans} one immediately has that, on $\ce^k$,
$$
\hat{\d}=\d  - (n-2k)\ii(\Up).
$$
 Thus this is conformally invariant if $k=n/2$ and also if
$k=0$ (since both sides then act trivially) but not otherwise. So
$\cH^0$, which is just the space $\cC^0$ of constant functions, is
conformally invariant and so also is $\cH^{n/2}$. Otherwise the
harmonics $\cH^k$ move in $\ce^k$ depending on the choice of metric.
The problem is to find a replacement space which is both isomorphic to
$H^k(M)$ and  stable under conformal transformations.

Without drawing on the details Hodge theory we can see at the outset
that the Riemannian system $(d,\d)$ has a finite dimensional null
space, since it is an elliptically coercive system.  The notion of
(graded) ellipticity we are using here, and below, is that the
operator concerned is a right factor of an operator with leading term
a power of the Laplacian. In this case the power is one:
$$
\horisys{\d }{d}\vertsys{d}{\d}=
-\Delta+\LOT,
$$ where $\Delta$ denotes the Bochner Laplacian $\nd^a\nd_a$ and
$\LOT$ indicates lower order terms.  This analysis suggests another
problem.

\noindent{\bf Problem 2: Gauge companion operators.} For each $k$ attempt 
to find a differential operator $G_k$ satisfying the following: \\ 
1. $G_k$ is conformally invariant on the null space
$\cC^k:=\cN(d:\ce^k\to \ce^{k+1})$. \\
2. (In a choice of scale) the system $(d,G_k)$ is elliptic. \\
 
On even dimensional manifolds there is another family of conformally
invariant operators between differential forms. These are operators:
$$
L_k:\ce^k\to \ce_k \quad L_k= (\d d)^{n/2-k} +\LOT \quad k\in
\{0,1,\cdots,n/2\}
$$
which we will term the {\em long} operators since they complete the
de Rham complex to the picture in figure \ref{figure}.  
While for
$k\geq 1$ the existence of operators between the bundles concerned can be concluded from the
general results of Eastwood and Slov\'ak in \cite{EastSlo}, these are not
unique and it turns out that our problems above are related to the
possibility of a special class of such operators.

\begin{figure}
\caption{\label{figure}
The conformal de Rham diagram in even dimensions}
 \begin{center}\begin{picture}(436,30)(0,-22)
\put(10,0)
{$\ce^0\stackrel{d}{\to}\cdots \stackrel{d}{\to}
\ce^{n/2-2}\stackrel{d}{\to}\ce^{n/2-1}\stackrel{d}{\to}
\ce^{n/2}\stackrel{\star
d}{\to}\ce_{n/2-1}\stackrel{\delta}{\to}\ce_{n/2-2}
\stackrel{\delta}{\to} \cdots\stackrel{\delta}{\to}\ce_0,
$}                                                         
\put(120,-14){\line(1,0){86}}
\put(120,-14){\line(0,1){7}}                               
\put(206,-14){\vector(0,1){7}}
\put(148,-10){\scriptsize$L_{n/2-1}=\d d$}

\put(71,-23){\line(1,0){184}}
\put(71,-23){\line(0,1){16}}                               
\put(255,-23){\vector(0,1){16}}
\put(93,-19){\scriptsize$L_{n/2-2}$}

\put(12,-31){\line(1,0){320}}
\put(12,-31){\line(0,1){24}}                               
\put(332,-31){\vector(0,1){24}}
\put(44,-27){\scriptsize$L_{0}$}

\end{picture}\end{center}
\end{figure} 

\noindent{\bf Problem 3: A preferred class of long operators.} The problem 
here is to attempt to establish existence of, or even better give a 
construction for, conformally invariant operators $L_k:\ce^k\to \ce_k$
which factor through $\d$ and $d$ in the sense,
\begin{equation}\label{factform}
L_k= \d M_{k+1} d .
\end{equation}
Such operators exist in the conformally flat case and we want the $L_k$ to
generalise these flat case operators. The operators sought 
should also be natural, that is in a choice of scale they should be given 
by a formula  polynomial in the Levi Civita connection 
$\nabla$ and its curvature $R$.

For $k=n/2-1$ the solution to this problem is the Maxwell operator $\d
d$ mentioned above.    Somewhat more compelling
evidence that there could be a positive solution to this problem, in general,
 dates
back to Branson's \cite{tbms} which provides direct constructions of
such operators at orders 4 and 6 (i.e. on, respectively, $\ce^{n/2-2}$
and $\ce^{n/2-3}$).  At the other extreme of order Graham et al (GJMS)
\cite{GJMS} give an order $n$ conformally invariant differential
operator $P_{n}:\ce^0\to \ce_0$ which has the desired form
\nn{factform}.

A powerful, and well understood, tool for generating conformally
invariant operators from other appropriate conformally invariant
operators is the {\em curved translation principle} of Eastwood and
Rice \cite{EastRice,Esrni}. However this does not predict the
operators $P_{n}$. The existence of these is subtle and the
construction of GJMS uses the ambient metric of Fefferman-Graham. We
will see below that the curved translation principle cannot in general
yield operators of the form \nn{factform}. Before we discuss these
difficulties, we draw in the final main problem which is linked to the
GJMS operators.

 For each even integer $n_0\geq 4$, the construction of GJMS gives,
not just an operator in dimension $n_0$, but an order $n_0$
conformally invariant differential operator
$P_{n_0}=\Delta^{n_0/2}+\LOT$ for all odd dimensions and for even
dimensions $n\geq n_0$. It is an observation of Branson
\cite{tbprsnl,tomsharp} that, in any choice of scale, these take the
form $\d M_1 d+\frac{1}{2}(n-n_0)Q(n,n_0)$, also that there exists a
universal expression for the order zero part $Q(n,n_0)$ which is
rational in $n$ (without singularity at $n_0$) and that setting
$n=n_0$ in this yields a remarkable curvature quantity
$Q:=Q(n_0,n_0)$. This has the conformal transformation
\begin{quote}
(I) $\quad Q^{\hat g}=Q^g+ P_{n_0} \om   \quad \mbox{where} \quad \hat g : = e^{2\om} g .$
\end{quote}
In dimension $n_0$, $P_{n_0} = \d M_1 d $
and so 
$
\cR(P_{n_0})\subseteq \cR(\d)
$. 
Thus $Q$ gives a conformally invariant operator
\begin{quote}
(II) $\quad Q: \cC^0\to  H_0(M)\cong H^n(M),$
\end{quote}
by $c\mapsto [c Q] $.  Also $Q$ has density weight $-n_0$ and so, since 
$ \cR(P_{n_0})\subseteq \cR(\d) $
and $Q$ transforms conformally as in (I), it follows that $\int Q$ is a conformal invariant. In the
conformally flat case $\int Q$ recovers a non-zero multiple of the
Euler characteristic and so the map (II) is in general non-trivial.
In fact the operators $P_{n_0}$ are (formally) self-adjoint so more
generally we have,
\begin{quote}
(III)  $\quad  c\in \cC^0 \mbox{ and } u\in \cN(P_{n_0}) ~\Rightarrow ~ 
\int u Q c
~\mbox{ is conformally invariant.}
$
\end{quote}

It turns out that the Q-curvature has a serious role in geometric
analysis and low dimensional topology \cite{CQY,CY}. There are also
connections with the AdS/CFT correspondence of quantum gravity
\cite{FGrQ} and scattering theory \cite{GrZ}. There have been recent
alternative direct constructions of $Q$ via tractor calculus
\cite{GoPet}, and the ambient construction \cite{FeffH}, which avoid
dimensional continuation. However there are still many mysteries. In
particular an important question is  whether there are other similar quantities.

\noindent{\bf Problem 4: Understand/generalise Branson's Q-curvature.}
Broadly the problem here is to find an analogue or ``closest
relative'' of the Q-curvature for forms. Aside from shedding light on
the Q-curvature itself, the idea is to investigate the existence of
other objects which are not conformally invariant locally and yet, by
some analogy with (I), (II) and (III) above, yield new global 
conformal invariants. 

\subsection{Earlier work} The 4$^{\rm th}$ order GJMS operator $P_2$ is first due to 
Paneitz \cite{Pan}. In dimension 4 this acts on functions and (given a
choice of metric $g$) has the formula
$$
\d(d\d+2\J-4\V\sharp)d .
$$ Here $\V$ is the Schouten (or Rho) tensor, viewed as weighted
section of End$(T^*M)$, $\J$ its trace and $\sharp$ is the obvious
tensorial action.  Recall these are related to the Ricci tensor for
$g$ by $\Ric (g)=(n-2)\V+\J g$.  Thus, in dimension 4, the operator
$G=\frac12\d(d\d+2\J-4\V\sharp) $ is conformally invariant on exact
1-forms.  It is easily verified by direct calculation that $G$ is also
conformally invariant on the space of closed 1-forms and hence also on
the null space of the Maxwell operator. Eastwood and Singer made this
observation and proposed $G$ as a gauge operator for the Maxwell
operator \cite{EastSin1}.  Note that
$$
\horisys{\d d\d}{d}\vertsys{ d}{\d d\d+\LOT}=
\Delta^2+\LOT.
$$
and so $(d,~G)$ gives a solution to problem 2 for $d$ on 1-forms in
dimension 4.

It was observed by the author and Branson \cite{opava} that this gauge
operator can be recovered from an adaption of the curved translation
principle. This idea also provides a conceptual framework and
practical approach to constructing gauge operators for many other
conformally invariant operators. First we need the basic idea of the 
tractor connection
and its associated calculus.

Any conformal $n$-manifold, such that $n\geq 3$, admits a unique
normal tractor bundle and connection.  The tractor connection is a connection 
on a vector bundle that we term the standard conformal tractor bundle
$\bT$. We write $\ct$ for the space of sections of $\bT$.  For a
choice of metric $g$ from the conformal class this bundle can be
identified with the direct sum $ [\bT]_g= E[1]\oplus E^1[1] \oplus
E[-1]$, where $E^1[1]$ means $T^*M\otimes E[1]$. Assigning abstract
indices we could instead write $ [\bT^A]_g= E[1]\oplus E_a[1] \oplus
E[-1]$.  Thus a section $V\in\ct$ then corresponds to a triple
$(\alpha,\mu,\rho)$ of sections from the direct sum according to $V^A=
Y^A\alpha+Z^{Ab}\mu_b+X^A\rho $ (where this defines the ``projectors''
$X$, $Z$ and $Y$).
Under a conformal rescaling $g\mapsto \widehat{g}=e^{2\w}g$, this triple transforms according to 
\begin{equation}\label{matrixtr}
[V]_g=\left(\begin{array}{c}  \alpha \\
                              \mu_b \\
                              \rho   \end{array}
\right) \mapsto \left(\begin{array}{ccc}  1 & 0 & 0  \\
                              \Up_b & \delta^a_b & 0 \\
                              -\frac{1}{2}\Up^a\Up_a & -\Up^a & 1  \end{array}
\right) 
\left(\begin{array}{c}  \alpha \\
                              \mu_a \\
                              \rho   \end{array}
\right) =[V]_{\hat g}
\end{equation}
where $\Up:=d\om$. It is easily verified that this determines an
equivalence relation on the triples over the equivalence relation on
metrics and hence the quotient gives $\bT$ as a well defined vector
bundle on $(M,[g])$ with a composition series $ \bT= E[1]\lpl E^1[1]
\lpl E[-1]$ (meaning that $E[-1]$ is a subbundle of $\bT$ and $E^1[1]$
is a subbundle of the quotient $\bT/E[-1] $).

In terms of this splitting for $g$ the conformally invariant tractor
metric is given by $h(V,V)=\bg^{ab}\mu_a\mu_b +2 \alpha \rho$.  The
tractor connection \cite{BEGo} is given by
\renewcommand{\arraystretch}{1}
\begin{equation}\label{stdtracconn}
[\nabla_a V^B]_g =
\left(\begin{array}{c} \nabla_a \alpha-\mu_a \\
                       \nabla_a \mu_b+ \bg_{ab} \rho +\V_{ab}\alpha \\
                       \nabla_a \rho - \V_{ab}\mu^b  \end{array}\right) . 
\end{equation}
\renewcommand{\arraystretch}{1.5} In terms of this formula, the
conformal invariance of the connection is recognised by the fact that
the components on the right-hand-side transform, under conformal
rescaling, according to \nn{matrixtr}.  In subsequent calculations we
will often omit the $[\cdot]_g$ which emphasises the choice of
splitting, since, in any case, this should be clear by the context.  For
the purposes of calculations it is often more convenient (see
\cite{GoPet}) to use that the connection is determined by
\begin{equation}\label{connids}
\begin{array}{rcl}
\nd_aX_A=Z_{Aa}\,, &
\nd_aZ_{Ab}=-\V_{ab}X_A-Y_A\bg_{ab}\,, & \nd_aY_A=\V_{ab}Z_A{}^b ,
\end{array}
\end{equation}
and the Leibniz rule. The tractor bundle and connection are induced by, and are
equivalent to, the normal conformal Cartan connection, see \cite{CapGotrans}.

The bundle of $k$-form tractors $\bT^k$ is the $k^{\underline{\rm
th}}$ exterior power of the bundle of standard tractors. This has a
composition series which, in terms of section spaces, is given by
\begin{equation} \label{formtractorcomp}
\cT^k=\Lambda^k\cT\cong\ce^{k-1}[k]\lpl\big(\ce^k[k]\oplus
\ce^{k-2}[k-2]\big) \lpl\ce^{k-1}[k-2].
\end{equation}
Given a choice of
metric $g$ from the conformal class there is a splitting 
of this composition series 
corresponding to the splitting of $\bT$ as mentioned above.   
Relative to this, a typical $k$-form tractor field $F$
corresponds to a 4-tuple $(\alpha,\m,\f,\r)$ of sections of the direct sum 
(obtained 
by replacing each $\lpl$ with $\oplus$ in \nn{formtractorcomp}) and we write
$$
F=\bbY^k\cdot \alpha +\bbZ^k \cdot \m + \bbW^k \cdot \f+ \bbX^k\cdot \r,
$$
where `${\cdot}$' is the usual pointwise form inner product in the
tensor arguments,
$$
\f\cdot\psi=\dfrac1{p!}\f^{a_1\cdots a_p}\psi_{a_1\cdots a_p}\ 
\mbox{ for $p$-forms},
$$
and for $k>1$, if $\uw$ denotes the wedge product in the tractor arguments, then we have
\begin{equation}\label{parts}
\bbZ^k=Z\uw\bbZ^{k-1},\ \ \bbX^k=X\uw\bbZ^{k-1},\ \
\bbY^k=Y\uw\bbZ^{k-1},\ \
\bbW^k=Y\uw X\uw\bbZ^{k-2}.
\end{equation}
By convention, $\bbZ^0=1$ and $\bbZ^{-1}=0$. 
The connection on $\bT$
gives a (conformally invariant) connection on $\bT^k$ by the Leibniz rule.
Under a change of scale $\hat g=e^{2\om}g$, it follows from the 
transformation law for the standard tractor bundle \nn{matrixtr}, that
\renewcommand{\arraystretch}{1.1}
$$
\begin{array}{rl}
\widehat \bbX&=\bbX, \\
\widehat \bbZ&=\bbZ+\e(\Up)\bbX, \\
\widehat \bbW&=\bbW-\ii(\Up)\bbX, \\
\widehat\bbY&=\bbY-\ii(\Up)\bbZ-\e(\Up)\bbW
+\frac12(\e(\Up)\ii(\Up)-\ii(\Up)\e(\Up))\bbX,
\end{array}
$$ 
\renewcommand{\arraystretch}{1.1} where again $\Up = d\omega$ and
the interior and exterior multiplication apply to the tensor indices
of the $\bbX$, $\bbZ$, $ \bbW$ and $\bbY$ projectors.

We are now ready to investigate gauge operators for the Maxwell
operator on $\ce^{n/2-1}$. By \nn{formtractorcomp} the forms $\ce^k$
turn up at the $\bbZ^{k}$ slot of $\ct^{k}[-k]$. It is straightforward
to verify, using the formulae above, that
\begin{equation}\label{splits}
\mu \mapsto 
\left(\begin{array}{c} 0 \\ \m \\ (n-2)^{-1}\d\m
\end{array}\right), ~ k=1 \hspace*{2mm} \rm{and} \hspace*{2mm}
\mu \mapsto 
\left(\begin{array}{c} 0 \\ \m\qquad\ 0\ \\ (n-2k)^{-1}\d\m
\end{array}\right), ~ k\geq 2
\end{equation}
 are conformally invariant differential splitting operators $S_k:
\ce^{k}\to \ct^{k}[-k]$.

Now recall the Yamabe operator (or conformal Laplacian) 
$$\Box =
-\nabla^a\nabla_a -(1-n/2) \J
$$ is conformally invariant on the densities $\ce[1-n/2]$. In fact
this same formula also gives an invariant operator on the sections of
$E[1-n/2]\otimes \bU $, where $\bU$ is any vector bundle with
connection, provided we view $\nabla$ as the coupled Levi Civita
vector bundle connection. Operators with this property are said to be
{\em strongly invariant}. (This is just a slight variation of the
notion introduced in \cite{Esrni}.) To see this one can simply observe
that the direct calculation, using \nn{densconntrans} and
\nn{conntrans}, which verifies the invariance of $\Box$ does not
require the commutation of covariant derivatives.
In particular we may couple with the tractor bundle. As an
immediate application, note that $S_{n/2-1}$ takes values in
$\ct^{n/2-1}[1-n/2]$ and so the composition $\Box
S_{n/2-1}:\ce^{n/2-1}\to \ct^{n/2-1}[-1-n/2]$ is conformally
invariant.

In dimension 4, for example, $\mu$ is a 1-form and we should
(following \cite{opava}) calculate $(-\nabla^a\nabla_a + \J) (Z\cdot
\mu + X\frac12 \d \mu)$ using \nn{connids}.  This yields
$$
\left(\begin{array}{c}0\\
\d d \m \\
\frac12\d(d\d\m+2\J-4\V\sharp)\m
\end{array}\right) = \left(\begin{array}{c}0\\
\d d \m \\
G \m
\end{array}\right). 
$$ We have recovered exactly the Maxwell operator and Eastwood-Singer
gauge pair.  Since this construction is conformally invariant it is
immediate from this final formula that $G$ is invariant on the null
space of the Maxwell operator $\d d$.  This construction is in the
spirit of the {\em curved translation principle} of Eastwood et al.\
\cite{EastRice,Esrni}; we have obtained the (operator,~gauge) pair by
translating from the Yamabe operator.

Buoyed by this success we are drawn to immediately try the same idea
to obtain a gauge operator for the Maxwell operator in higher even
dimensions. However this time we find that, for $\mu\in \ce^{n/2-1}$, we get 
$$
\Box \left(\begin{array}{c} 0 \\ \m\qquad\ 0\ \\ (n-2k)^{-1}\d\m
\end{array}\right) = 
\left(\begin{array}{c} 0\\
(\d d+C\sharp\sharp)\m \qquad 0 \\
\frac12\d d\d\m+\d(\J\m)-2\V\sharp\d\m-2\d(\V \sharp \m)+\frac12 C\sharp\sharp\d\m
\end{array}\right),
$$ where $C$ is the Weyl curvature of the conformal structure (recall
the Weyl curvature is conformally invariant) and in the
$C\sharp\sharp$ action we view this as a (weighted) section of the tensor square
of $\End(T^*M)$.  We see here that, unfortunately, the Maxwell operator does not
automatically turn up as the leading slot. We did not encounter this
in dimension 4 because when $\mu$ is a 1-form $C\sharp\sharp \mu$
vanishes (since $C$ is trace-free). Of course $ C\sharp\sharp\m$ is
conformally invariant but we cannot simply subtract this and maintain
conformal invariance without also adjusting the $\bbX$-slot.  The
existence of a correction, in such circumstances, is a delicate
matter. It turns out that in this instance there is a fix. The output
above decomposes into a sum of conformally invariant tractors
according to
\begin{equation}\label{splitmax}
\left(\begin{array}{c}0\\
\d d\m \qquad 0 \\
\frac12\d d\d\m+\d(\J\m)-2\d(\V\sharp\m)
\end{array}\right)+
\left(\begin{array}{c}0\\
C\sharp\sharp\m \qquad 0 \\
Y\bullet\m
+\frac12C\sharp\sharp\d\m
\end{array}\right),
\end{equation}
where $Y=Y^{ab}{}_c:=\nd^a\V^b{}_c-\nd^b \V^a{}_c$ and $Y\bullet\m$
means $-\sum_{s=2}^k Y^{a_1 b}{}_{a_s}\m_{a_1\cdots b\cdots a_k}$.

So $G_{n/2-1}:=\frac12\d d\d+\d(\J)-2\d(\V\sharp)$ is a gauge
operator for the Maxwell operator (and so also for $d$) in all even
dimensions. By construction it is conformally invariant on
$\cN(d)\subseteq \cN(\d d)$ (in fact this is an equality in the
compact Riemannian setting) and from our earlier observations it
combines with $d$ to give an elliptic system. Finally note that if we
apply the Maxwell operator plus gauge system to an exact $\mu =d \nu$ 
we obtain 
$$
\left(\begin{array}{c}0\\
0 \qquad 0 \\
G_{n/2-1}d\nu
\end{array}\right), 
$$ 
since $\d d$ annihilates exact forms.  In the alternative notation this is
$$\bbX\cdot G_{n/2-1}d\nu= \bbX\cdot \big( \frac12  \d
d\d\m+\d(\J\m)-2\d(\V\sharp\m)d\nu\big) .
$$ Since this is conformally invariant by construction,
$\widehat\bbX=\bbX$, and the coefficients in the other slots are all
zero, it follows that
\begin{equation}\label{L4}
\big(G_{n/2-1}d = \frac12\d(
d\d +2\J-4\V\sharp)d\big): \ce^{n/2-2}\to \ce_{n/2-2}
\end{equation} 
is a conformally invariant long operator of the form proposed in
problem 3. 

Although we have succeeded in pushing this calculation through there
are are two main problems which suggest that this approach would be
difficult, if not impossible, to generalise sufficiently to deal with
our problems 2 and 3 in general.  One is that even at the low order of
example treated, the calculations leading to the results presented were
non-trivial and involved, for example, the Bianchi identity
$\nd_{[a}R_{bc]de}=0$. More seriously the decomposition in
\nn{splitmax} involved inspecting the explicit formulae and solving
equations to extend the $C\sharp\sharp\mu$ term to a conformally
invariant tractor operator. Ab initio one does not know that this will
succeed.  The weight of $C\sharp\sharp\mu$ is exactly such that the
standard tools using Lie algebra cohomology \cite{CSSannals} or
central character as in the theory of Verma modules \cite{EastSlo} fail
to indicate the existence of this extension.

As a final point, in this lecture, let us note that the operator
\nn{L4}, and the proposed higher order analogues, are \underline{not}
strongly invariant.  We can easily see this directly for \nn{L4}. Let
us write $G:= G_{n/2-1}$ and note first that, from the transformation
formula $\hat \bbZ=\bbZ+\e(\Up)\bbX$ and the invariant splitting
\nn{splitmax}, it follows that $\widehat{G}=G-\ii(\Up)\d d$.  This is
also readily verified by direct calculation using the conformal
transformation formulae \nn{densconntrans} and \nn{conntrans} for the
Levi Civita connection.  Now consider a coupled variant of $G$ acting
on a vector bundle valued $(n/2-1)$-form $\mu$. Suppose that the
vector bundle has a connection $A$, with curvature $F$, and $G^A$ is
given by the formula above for $G$, except that $d$ and $\d$ are
replaced by their connection coupled variants $d^A$ and $\d^A$.  Now
the direct computation of the conformal transform of this,
$\widehat{G}^{A}\mu$, is the same as for the case of forms except that
now vector bundle curvature terms may enter from the commutation of
derivatives. Given that $G^A$ is just a 3${}^{\rm rd}$ order operator
one easily sees that
$$
\widehat{G}^A\mu =G^A\mu-\ii(\Up)\d^A d^A \mu + \Up\cdot F\cdot \mu,
$$ where $\Up\cdot F\cdot \mu$ indicates a sum of terms linear in
$\Up$, $ F$ and $ \mu $.  Now let us suppose
that $\mu =d^A\nu$ where the vector bundle valued $(n/2-2)$-form
$\nu$ satisfies $\nabla^A \nu (p)=0$ for some point $p\in M$. Note that
$d^A$ is conformally invariant on $\nu$. Then, at $p$, $\mu$
vanishes and we have
$$ \widehat{G}^Ad^A \nu =G^Ad^A\nu-\ii(\Up)\d^A
d^Ad^A \nu= G^Ad^A \nu -\ii(\Up)\d^A F\wedge \nu ,
$$ where $F\wedge \nu$ includes an implicit curvature action on $\nu$.
It is an elementary matter to verify by example that the term
$\ii(\Up)\d^A F\wedge \nu$ does not vanish in general. Thus $G^Ad^A$ is
not conformally invariant, and so \nn{L4} is not strongly invariant.

The operator \nn{L4} is not unique. For example we could add to it the
conformally invariant term $\d C\sharp\sharp d$. It is natural to
wonder if there is some modification which still has the form $\d M d$
but which is strongly invariant. One needs to be careful considering
such arguments since strong invariance is really a property of the
{\em formulae} for operators rather the operators themselves.
Nevertheless we will show that in fact, apart from the $2^{\rm
  nd}$-order Maxwell operators $\d d$, none of the operators sought in
problem 3 can be strongly invariant. (That is there are not strongly
invariant formulae for these operators which have the form $\d M d$.)
This means that they cannot be obtained by the usual use of the curved
translation principle since that procedure involves composing strongly
invariant operators to obtain new strongly invariant operators.  This
is an important feature of the desired operators, so we state the
result as a proposition. (In fact we give a stronger result.) In
proving this we will use, what is now a well known result (which can
be deduced from the results in \cite{EastSlo}) as follows.  On the
conformally flat sphere one has the invariant operators $\ce^{n/2-1}
\stackrel{\star d}{\to} \ce^{n/2} $ and $\ce^{n/2}\stackrel{\d}{\to}
\ce_{n/2-1}$. These with the operators indicated in figure
\ref{figure}, give, up to linear combinations, all of the conformally
invariant operators between the differential form bundles in figure
\ref{figure}. It follows, for example, that the composition of
operators from the figure always yields a trivial operator.
\begin{proposition}\label{notstrong}
Suppose that for $k\in \{0,1,\cdots ,n/2-2\}$
$$
 S d :\ce^{k}\to \ce_k
$$
is a strongly invariant natural conformally invariant operator.  Then 
this operator vanishes on conformally flat structures. 
\end{proposition}
\noindent{\bf Proof:}
Suppose first that $k\geq 1$. Let $\cV$ be a trivial bundle, with
fibre $V$. Let us equip this with a family of connections which differ
from the trivial connection by $tA$, where $t$ is a real parameter and
$A$ is any field of $\End(V)$-valued 1-forms.  Since $ S d$ is
strongly invariant we can couple to the connection corresponding to
$tA$, for each $t$, to obtain the conformally invariant operator $
S^{tA} d^{tA}$ on $\cV$-valued $k$-forms. Since also the exterior
derivative is strongly invariant it follows that the composition $
S^{tA} d^{tA} d^{tA}$ is also conformally invariant on $\cV$-valued
$(k-1)$-forms. But this is a non-zero multiple of 
$$
 S^{tA} F^{tA}
$$
where $F^{tA}$ is the curvature of the connection $tA$. 
Of course the 
$F^{tA}$ acts by the exterior product, via its form indices, as well
as the usual adjoint action of a curvature.  Viewing $t$ as a
parameter, it is clear that the displayed operator can be expressed by
a formula polynomial in $t$ and so its derivatives, with respect to
$t$, are also conformally invariant. In particular if we write $F_0:=d
F^{tA}/dt|_{t=0}$ then, by differentiating and evaluating at 0,
we obtain that
$$
S F_{0}
$$
is conformally invariant on $\cV$-valued $(k-1)$-forms.  Now
$F_{0}$ is an $\End(\cV)$-valued 2-form. Since $\End(V)$ is
canonically isomorphic to its dual, we may view $F_{0}$ instead as map
from $\End(\cV)\to \ce^2$. It is easily verified that, by suitable
choice of $V$ and the field $A$, one can arrange that this map is
surjective. So let us assume this. We have stated that if $H$ is any
$\cV$-valued $(k-1)$-form then $S F_{0} H$ is conformally invariant.
Now suppose $W$ is a section of $\cV^*$ that is parallel for the trivial
connection on $\cV^*$. Then $W\cdot S F_{0} H = S (F_{0} H)\cdot W$, where the 
`\cdot' indicates that the section $W$ is contracted into the free 
$\cV$-index of 
$ (F_{0} H)$. Thus $S$ is conformally invariant on the $(k+1)$-form 
$(F_{0} H)\cdot W $. Since $S$ is linear, it is also 
conformally invariant on sums of $(k+1)$-forms constructed this way and so, 
by the surjectivity of $F_{0}$, we can conclude that  
$ S:\ce^{k+1}\to \ce_{k}$ is conformally invariant. But
then it follows that this is trivial in the flat case because, for
$k$ in the range assumed (from the classification described above),
there are no non-trivial conformally invariant operators, on the
conformal sphere, of the form $ \ce^{k+1} \to \ce_{k}$.  This does the
cases $k\neq 0$.

Now suppose, with a view to contradiction that $ L: \ce^0\to
\ce_0$ is natural, strongly invariant, and is non-trivial on the conformal
sphere. Then the  leading term is $\Delta^{n/2}$, at least up to a constant non-zero multiple. Coupling to the standard
tractor bundle and connection we may conclude the
existence of a conformally invariant operator
$$
H^B{}_A:\ct^A\to \ct^B
$$
with principal part $\Delta^{n/2}$ (where now $\Delta$ is the
tractor-coupled Laplacian). Now there is the so-called tractor-D
operator \cite{BEGo} which is (strongly) conformally invariant and  
given by the 
formula 
$$
D^A f:=(n+2w-2)w Y^A f+ (n+2w-2)Z^{Aa}\nabla_a f -X^A(\Delta+w\J) f
$$
for $f$ any weight $w$ tractor or density field.  
Composing first with this on the right
it is easily verified (or see \cite{gosrni} or \cite{GoPet}) that, in the conformally flat case, we have 
$$
H^B{}_A D^B f=-X^B \Delta^{n/2+1} f
$$
for $f$ any weight $1$ density field. From this it follows easily that, 
in the general curved setting, the conformally invariant composition
$$
D_B H^B{}_A D^A:\ce[1]\to \ce[-1-n] .
$$
has leading term a
non-zero constant multiple of $\Delta^{n/2+1}$.  However this is a
contradiction as there is no such operator \cite{GHnon}.  \quad $\Box$

\noindent Notice that we have proved a little more than what is 
claimed in the proposition. We have shown that there is no strongly
conformally invariant curved analogue of the operator $P_n :\ce^0\to
\ce_0$ on the sphere, regardless of its form.

\section{Lecture 2 -- Operators like Q and the ambient connection}

In the first half of this lecture we show that the 4 problems are
 solved simultaneously (with some mild qualifications) by a sequence
 of remarkable operators which include and generalise the
 $Q$-curvature.  The construction of these operators uses the
 Fefferman-Graham ambient metric construction and its relationship to
 the tractor calculus. In the second half of the lecture we set up the
 background for this.

\subsection{The solution}

We collect the main points into a theorem which
includes some of the central results in \cite{deRham}:
\begin{theorem}\label{Qthm} In each even dimension $n$ there exist 
natural Riemannian differential operators 
$$
Q_k:\ce^k\to \ce_k, \quad \mbox{non-zero for} \quad k=0,1,\cdots , n/2, 
$$ 
(and for other $k$  we take these to be zero)
with the following properties:\\
\noindent  \IT{i}
$ Q_0 1$ is the Branson Q-curvature.\\
\noindent  \IT{ii}  As an operator on closed $k$-forms
$Q_k $ has the conformal transformation law
$$
\widehat{Q}_k=Q_k  + \d Q_{k+1} d \om
$$ 
where $ \hat{g}=e^{2\om} g$, for a smooth function $ \om$, and on the
right-hand-side we view $\om$ as a multiplication operator. \\
\noindent \IT{iii} 
$$
\d Q_{k+1} d: \ce^{k} \to \ce_{k} \quad k=0,1,\cdots , n/2 -1
$$ 
is conformally invariant (from \IT{ii}), formally self-adjoint and has leading term a
non-zero multiple of $(\d d)^{n/2-k}$. \newline
\noindent  \IT{iv} The system $(\d Q_{k+1} d,~ \d Q_k)$ is elliptic 
($k\in \{0,1,\cdots , n/2-1\}$) and, for each $k$, 
$\d Q_k$ has the conformal transformation 
$$
\widehat{\d Q_k}= \d Q_k + c\ii(d\om)  \d Q_{k+1} d, 
$$
where $\hat g= e^{2\om}g$ and $c$ is 
a constant.
In particular 
$
\d Q_k
$
is conformally invariant  on $\cN(\d Q_{k+1} d)$.
\end{theorem}

By part \IT{iii}, our problem 3 is solved by taking $L_k= \d Q_{k+1}
d$. 
It follows immediately that in each dimension there is a family of conformally invariant 
{\em detour complexes},
\begin{equation}\label{detour}
\ce^0\stackrel{d}{\to}\cdots\stackrel{d}{\to}\ce^{k-1}\stackrel{d}{\to}
\ce^k\stackrel{L_k}{\to}\ce_k\stackrel{\delta}{\to}\ce_{k-1}
\stackrel{\delta}{\to}
\cdots\stackrel{\delta}{\to}\ce_0,
\end{equation}
which generalise the Maxwell detour complex.
Since $L_k$ has leading term $(\d d)^{n/2-k}$ these are elliptic
(i.e.\ exact at the symbol level).

Clearly $\cN(d) \subseteq \cN(\d Q_{k+1} d)$ and so $\d Q_k$ is
conformally invariant on $\cN(d)$. 
Since $\d Q_{k+1} d$ has leading
term $(\d d)^{n/2-k}$, it follows that the ellipticity of the system
$(\d Q_{k+1} d,~ \d Q_k)$, as asserted in \IT{iv}, implies that $(d,\d
Q_k)$ is elliptic. So setting $G_k :=\d Q_k$ gives a solution to
problem 2. Thus, writing $\cC^k$ for the space of closed $k$-forms, we propose 
$$
\cH^k_G :=\cN(G_k:\cC^k \to \ce_{k-1})  
$$ as the space of {\em conformal harmonics}, for $k=0,1,\cdots,
n/2$. Since $G_k$ is conformally invariant on $\cC^k$ it follows that
this space is conformally invariant.

Of course property \IT{ii} generalises the transformation
formula (I) of the Q-curvature (where we view $Q$ as a
multiplication operator on the constant functions). Then note that 
if $c\in \cC^k$, and  $u\in \ce^k$ then 
$$
\int_M (u, Q^{\hat g}_k c)d\mu_{\sbg} =\int_M (u, Q^{g}_k c+ L_k \om c) d\mu_{\sbg} 
= \int_M (u, Q^{ g}_k c) d\mu_{\sbg}+ \int_M (L_k u, c)d\mu_{\sbg}
$$
as $L_k$ is formally self-adjoint. Here we are using
$(\cdot,\cdot)$ for the complete contraction, via $\bg^{-1}$, of forms. 
So if $u\in\cN(L_k) $
then the last term vanishes and we have
$$ 
  c\in \cC^k \mbox{ and } u\in \cN(L_k) ~\Rightarrow ~ 
\int (u, Q_k c) d\mu_{\sbg}
~\mbox{ is conformally invariant.}
$$
which generalises property (III) of the $Q$-curvature.  

The
transformation law in \IT{i} implies that $Q_k$ gives a conformally
invariant map $Q_k:\cC^k \to \ce^k/\cR(\d)$. If $u\in \cH^k_G$,
 then $u$
is both closed and in the null space of $G_k$, and so $\d Q_k u =0$
since $G_k=\d Q_k$. Thus $Q_k$ fives a conformally invariant map
\begin{equation}\label{cohomap}
Q_k: \cH^k_G \to
H_k(M)=\cN(\d:\ce_k\to\ce_{k-1})/\cR(\d:\ce_{k+1}\to\ce_k)\cong H^k(M).
\end{equation}
Note that since $Q_0$ takes values in densities, $G_0=\d Q_0$ is
trivial and $\cH^0_G=\cC^0$ and so the result displayed generalises to
the $Q_k$ property (II) of the Q-curvature. Examples where the maps
\nn{cohomap} are non-trivial are given in \cite{deRham}.

It remains to check to how accurately the defined conformal harmonics
reflect the de Rham cohomology.  We have already observed that
$\cH^0_G=\cC^0\cong H^0(M)$. It turns out that $Q_{n/2}$ is a
non-vanishing constant (as a multiplication operator) and so at middle
forms we recover the usual harmonics, $\cH^{n/2}_G=\cH^{n/2}\cong
H^{n/2}(M)$.  Between these extremes, it
is easy to obtain an estimate on the size of the space $\cH^k_G$.
Note there is a map $\cH^k_G\to H^k(M)$ given by mapping closed forms
to their class in cohomology.  If $ w \in \cH^k_G$ is mapped to the
class of 0 in $H^k(M)$ then $w$ is exact. Say $w=d\f$. Since, in
addition, $\d Q_{k} w=0$, it follows that $\d Q_{k} d\f=0 $, that is
$L_{k-1}\f =0$.  Now $L_{k-1} = \d Q_{k} d$, so
$\cC^{k-1}\subseteq \cN(L_{k-1})$ and it follows that there is an exact sequence
$$
0\to \cC^{k-1}\to \cN(L_{k-1})  \stackrel{d}{\to} \cH^k_G \to H^k(M). 
$$
Now $H^{k-1}(M)$ is the image of
$\cC^{k-1} $ under the composition $\cC^{k-1}\to \cN(L_{k-1})\to
H^{k-1}_L(M) $, so finally we obtain
$$
0\to H^{k-1}(M)\to H^{k-1}_L(M) \to \cH^k_G \to H^k(M) 
\quad \mbox{ for } k=1,\cdots ,n/2-1, 
$$
and so 
$$
\dim \cH^k_G \leq b^k + \dim (H^{k-1}_L(M)/ H^{k-1}(M))
$$ 
where $\dim H^{k}_L(M)$ is the cohomology at $\ce^k$ of the
sequence \nn{detour} and $b^k$ is the $k^{\rm th}$-Betti number, i.e.\
$b^k= \dim H^k(M)$. Obtaining a lower bound is not so
straightforward. Nevertheless, using Hodge theory and the map
\nn{cohomap} it can be shown \cite{deRham} that  $b^k \leq\dim
\cH^k $. Thus to have $\dim \cH^k_G = b^k $ it is sufficient for the
conformal regularity condition $H^{k-1}_L(M) = H^{k-1}(M)$ (or
equivalently $\cN(L_{k-1})=\cC^{k-1}$) to be satisfied. For $n=4$ and
$k=1$ this is the notion of strong regularity proposed by
\cite{EastSin2}.  Although, for each $k$, the regularity should hold
generically, in some appropriate sense, for compact conformal
Riemannian manifolds there are counter-examples to strong regularity
on 4-manifolds \cite{Sin}. In \cite{deRham} it is shown that there is
a condition weaker than $H^{k-1}_L(M) = H^{k-1}(M)$ which is necessary
and sufficient for $\dim \cH^k_G = b^k $.

\subsection{The Fefferman-Graham ambient construction}

Recall that $\pi:\cq\to M$ is the conformal bundle of metrics.
Let us use $\rho $ to denote the ${\Bbb R}_+$ action on $ \cq$ given
by $\rho(s) (x,g_x)=(x,s^2g_x)$.  An {\em ambient manifold\/} is a
smooth $(n+2)$-manifold $\aM$ endowed with a free $\Bbb R_+$--action
$\rho$ and an $\Bbb R_+$--equivariant embedding
$i:\cq\to\aM$.  We write $\X\in\frak X(\aM)$ for the fundamental field
generating the $\Bbb R_+$--action, that is for $f\in C^\infty(\aM)$
and $ u\in \aM$ we have $\X f(u)=(d/dt)f(\rho(e^t)u)|_{t=0}$.
If $i:\cq\to\aM$ is an ambient manifold, then an {\em ambient
metric\/} is a pseudo--Riemannian metric $\h$ of signature $(n+1,1)$
on $\aM$ such that the following conditions hold:\\
\smallskip
\noindent
(i) The metric $\h$ is homogeneous of degree 2 with respect to the
$\Bbb R_+$--action, i.e.\ if $\Cal L_{\sX}$
denotes the Lie derivative by $\X$, then we have $\Cal L_{\sX}\h=2\h$.
(I.e.\ $\X$ is a homothetic vector field for $h$.)

\noindent
(ii) For $u=(x,g_x)\in \cq$ and $\xi,\eta\in T_u\cq$, we have
$\h(i_*\xi,i_*\eta)=g_x(\pi_*\xi,\pi_*\eta)$. 
\smallskip
\noindent
Henceforth we will identify $\cq$ with its image in $\aM$ and suppress
the embedding map $i$.
\smallskip

In
\cite{FGast} Fefferman and Graham treat the problem of constructing a formal power
series solution along $ \Cal Q$ for the Goursat problem of finding
an ambient
metric $ \h$ satisfying (i) and (ii) and the condition that it be
Ricci flat, i.e.\ Ric$(\h)=0$. From their results and some minor subsequent observations 
\cite{GoPet,GJMS}
we have the following: there is a formal 
solution for $ \h$ satisfying (i), (ii) and with  
\renewcommand{\arraystretch}{1}
$$
\noindent {\rm (iii)}  \quad 
\Ric(\h)=0\; \left\{\begin{array}{l} \mbox{to all orders if $n$ is odd,} \\
\\
                                 \mbox{up to the addition of terms vanishing} \\
                                 \mbox{to order $n/2-1$ if $ n$ is even,}
\end{array} 
\right.
$$ with $\Q:= \h(\X,\X)$ a defining function for $\cq$ and
$\h(\X,\cdot)=\frac{1}{2}d \Q$ to all orders in both dimension
parities.  We will use the term ambient metric to mean an ambient
manifold with metric satisfying all these conditions.  
 Note
that if $M$ is locally conformally flat then the flat ambient metric
is a (canonical) solution to the ambient metric problem. It is
straightforward to check \cite{deRham} that this is forced in odd
dimensions while in even dimensions this extends the solution.  When
discussing the conformally flat case we assume this
solution.

We should point out that Fefferman and Graham give uniqueness statements for their
metric, but we do not need these here.  The uniqueness of the operators
we will construct is a consequence of the fact that they can be
uniquely expressed in terms of the underlying conformal structure as
we shall explain later.

We write $ \nda $ for the ambient Levi-Civita connection determined
by $ \h$ and use upper case abstract indices $A,B,\cdots $ for tensors
on $ \aM$. For example, if $ v^B$ is a vector field on $\aM $, then
the ambient Riemann tensor will be denoted $\aR_{AB}{}^C{}_{D}$ and
defined by $ [\nda_A,\nda_B]v^C=\aR_{AB}{}^C{}_{D}v^D$. In this
notation the ambient metric is denoted $ \h_{AB}$ and with its inverse
this is used to raise and lower indices in the usual way. Most often
we will use an index free notation and will not distinguish tensors
related in this way. Thus for example we shall use $\X$ to mean
both the Euler vector field $\X^A$ and the 1-form $ \X_A=\h_{AB}\X^B$.

The condition $\LX\h=2\h$ is equivalent to
the statement that the symmetric part of $\nda\X$ is $\h$.  On the other
hand, since $\X$ is exact, $\nda\X$ is symmetric.  Thus
\renewcommand{\arraystretch}{1.5}
\begin{equation} \label{ndaX} 
\nda\X =\h,
\end{equation} 
which in turn implies 
\begin{equation} 
\label{XRrel} 
\X\hook\aR=0. 
\end{equation} 
Equalities without qualification, as here, indicate that the results
hold either to all orders or identically on the ambient manifold.

Let $ \cce(w)$ denote the space of functions on $ \aM$ which are
homogeneous of degree $ w\in {\Bbb R}$ with respect to the action $
\rho$.  Recall that densities in $ \ce[w]$ are equivalent to functions
in $ \cce(w)|_\cq$. More generally (weighted) tractor fields
correspond to the restriction (to $\cq$) of homogeneous tensor fields
on $\aM$. A tensor field $F$ on $ \aM$ is said to be {\em homogeneous
of degree} $w$ if $\rho(s)^* F= s^w F$, or equivalently $ \cL_{\miniX}
F=w F$. The relationship between the Fefferman-Graham ambient metric
construction and the tractor connection was established in
\cite{CapGoFG}. Following this treatment we will sketch how the conformal
tractor bundle, metric and connection are related to the ambient metric.

 On the ambient tangent bundle $T\aM$ we define an action of $\Bbb
R_+$ by $s\cdot \xi:=s^{-1}\rho(s)_\ast \xi$. The sections of $ T\aM$
which are fixed by this action are those which are homogeneous of
degree $ -1$. Let us denote by $ \act$ the space of such sections and
write $\act(w)$ for sections in $\act\otimes \cce(w)$, where the
$\otimes$ here indicates a tensor product over $\cce(0)$.  Along $
\cq$ the $\Bbb R_+$ action on $T\aM$ agrees with the ${\Bbb R}_+$
action on $\cq$, and so the quotient
$(T\aM|_\cq)/\Bbb R_+$, yields a rank $ n+2$ vector bundle $\tilde{\bT}$ over
$\cq/\Bbb R_+=M$. By
construction, sections of $p:\tilde{\bT} \to M$ are equivalent to sections
from $\act|_\cq$. We write $\tilde\ct$ to denote the space of such sections.

Since the ambient metric $\h$ is homogeneous of degree $2$
it follows  that for vector fields $\xi$ and $\eta$ on $\aM$
which are homogeneous of degree $-1$, the
function $\h(\xi,\eta)$ is homogeneous of degree $0$ 
and thus descends to a smooth function on
$M$. Hence $\h$ descends to a smooth bundle metric $h$ of signature
$(n+1,1)$ on $\tilde{\bT}$. 

Next we show that the space $\act$ has a filtration reflecting the geometry of
$\aM$.  First observe that for $\phi\in \cce(-1)$, $ \phi \X\in
\act$. Restricting to $\cq$ this determines a canonical inclusion $
E[-1] \hookrightarrow \tilde{\bT} $ with image denoted by $\bV$. 
Since $ \X$ generates the fibres of $\pi:\cq\to M$ the smooth
distinguished line subbundle $ \bV\subset \tilde{\bT}$ reflects the
inclusion of the vertical bundle in $ T\aM|_\cq$. 
We write $X$ for the canonical section in $\tilde{\ct}[1]$ giving this
inclusion.  We define $ \bF$ to be the orthogonal complement of $ \bV$
with respect to $ h$.  Since $\Q=\h(\X,\X)$ is a defining function for
$\cq $ it follows that $ X$ is null and so $\bV\subset \bF$. Clearly $
\bF$ is a smooth rank $ n+1$ subbundle of $ \tilde{\bT}$. Thus $ \tilde{\bT}/ \bF$ is
a line bundle and it is immediate from the definition of $ \bF$ that
there is a canonical isomorphism $E[1]\cong \tilde{\bT}/ \bF $ arising from
the map $\tilde{\bT}\to E[1]$ given by $V\mapsto h(X,V)$. Now recall
$2\h(\X,\cdot)=d \Q$, so the sections of $\act|_\cq$ which correspond to
sections of $\bF$ are exactly those that take values in $ T\cq\subset
T\aM|_\cq$. 
Finally we note that if
$\tilde \xi$ and $\tilde{\xi}'$ are two lifts to $ \cq$ of $\xi\in
{\frak X}(M)$ then they are sections of $T\cq$ which are homogeneous
of degree 0 and with difference $\tilde \xi - \tilde{\xi}' $ taking
values in the vertical subbundle.  Since $\pi: \cq \to M$ is a
submersion it follows immediately that $\bF[1]/\bV[1]\cong TM\cong
T^*M[2]$ 
(where recall by our conventions $\bF[1]$ means $\bF\otimes E[1]$ etc.\/). 
Tensoring this with $E[-1]$ and combining this observation 
with our earlier results we can summarise the filtration of $\tilde{\bT}$ by
the composition series
\begin{equation} \label{trcomp} 
\tilde{\bT}= E[1]\lpl T^*M[1] \lpl E[-1] .
\end{equation}

Next we show that the Levi-Civita connection $\nda$ of $\h$ determines
a linear connection on $ \tilde{\bT}$. Since $\nda$ preserves $\h$ it follows
easily that if $U\in \act(w)$ and $ V\in \act(w')$ then $ \nda_U V \in
\act(w+w'-1)$.  The connection $ \nda $ is torsion free so $ \NX U -
\nda_U \X-[\X,U]=0$ for any tangent vector field $ U$.  Now $\nda_U \X
=U$, so this simplifies to $ \NX U = [\X,U]+U$.  Thus if $ U\in \act$,
or equivalently $[\X, U] =-U$, then $\NX U =0$. The converse is clear
and it follows that sections of $ \act$ may be characterised as those
sections of $ T\aM$ which are covariantly parallel along the integral
curves of $ \X$ (which on $ \cq$ are exactly the fibres of $\pi$).
These two results imply that $ \nda $ determines a connection $ \nd$
on $\tilde{\bT}$. For $ U\in \ct$, let $ \tU $ be the corresponding section of
$ \act|_{\cq}$. Similarly a tangent vector field $ \xi$ on $M$ has a
lift to a field $ \tilde{\xi}\in \act(1)$, on $ \Cal Q$, which is
everywhere tangent to ${\Cal Q}$.  This is unique up to adding $f\X$,
where $ f\in \cce(0)$. We extend 
$ \tU$ and $\tilde{\xi} $ smoothly and homogeneously to
fields on $ \aM$.  Then we can form $ \nda_{\tilde{\xi}} \tU$; this is
clearly independent of the extensions. Since $ \NX \tU =0$, 
the section $ \nda_{\tilde{\xi}} \tU$ is also
independent of the choice of $ \tilde{\xi}$ as a lift of $ \xi$. Finally, 
$ \nda_{\tilde{\xi}} \tU$ is
a section of $ \act(0)$ and so determines a section $ \nd_\xi U $
of $ \tilde{\bT}$ which only depends on $ U$ and $ \xi$. It is easily verified
that this defines a covariant derivative on $ \tilde{\bT}$ which, by
construction, is compatible with the bundle metric $ h$.

The ambient metric is conformally invariant; no choice of metric from
the conformal class on $M$ is involved in solving the ambient metric
problem.  Thus the bundle, metric and connection $(\tilde{\bT},h,\nd)$
are by construction conformally invariant.  On the other hand the
ambient metric is not unique (there is some diffeomorphism freedom and,
even allowing for this, recall that in even dimensions the
construction is only determined by the underlying conformal manifold
to finite order). Nevertheless it is straightforward to verify that $
\nd$ satisfies the required non-degeneracy condition and curvature
normalisation condition \cite{CapGotrans} that show that the bundle
and connection pair $(\tilde{\bT},\nd)$, induced by $\h$, is a normal
standard (tractor bundle, connection) pair.  So although the ambient
metric is not unique the induced tractor bundle structure
$(\tilde{\bT},h,\nd)$ is equivalent to a normal Cartan connection, and
so is unique up to bundle isomorphisms preserving the filtration
structure of $\tilde{\bT}$, and preserving $h$ and $\nd$. Hence we may drop the tildes and
identify $\tilde{\bT}$  with $\bT$ and $\tilde\ct$ with $\ct$.

Since $\ct$ corresponds to the ambient space $\act|_\cq$ and $\ce[w]$
corresponds to $\cce(w)|_\cq$ it follows, by taking tensor powers,
that homogeneous ambient tensors along $\cq$ are equivalent to
weighted tractor fields in the corresponding tensor power of $\ct$.
In particular this is true for exterior powers. 
The subspace of $\Gamma(\wedge^k T^*\aM)$ consisting of ambient
$k$-forms $F$ satisfying $\NX F =w F$ for a given $w\in {\Bbb R}$ will
be denoted $\act^k(w)$. We say such forms are (homogeneous) of {\em
weight} $ w$. Then we have that each section  $V\in \ct^k[w]$,  is  equivalent
to a section  $\tilde{V}\in \act^k(w)|_\cq$.

\subsection{Exterior calculus on the ambient manifold}

We need to identify which operators on the ambient manifold correspond
to, or determine, conformal differential operators on $M$. In
particular for our problems it turns out that operators on ambient
differential forms have a primary role. 

We will use $\ii(\cdot)$ and $\e(\cdot)$ as the notation for interior
and exterior multiplication by 1-forms on ambient forms, i.e.\ the
same notation as on $M$ and with the same conventions.  Thus for
example on differential forms, the Lie derivative with respect to $\X$
is given by $\LX=\ii(\X)\ad+\ad\ii(\X)$ and so its formal adjoint is
$\LX^*=\da\e(\X)+\e(\X)\da$.  Note that $\Q:=\h(\X,\X)$ may be alternatively expressed
$$
\Q=\ii(\X)\e(\X)+\e(\X)\ii(\X).
$$ 
It is useful for our calculations to extend the notation for
interior and exterior multiplication, in an obvious way, to operators
which increase the rank by one. For example, writing $\ad$ and $\da$
for respectively the ambient exterior and its formal adjoint, we have
$\ad\f=\e(\nda)\f$ and $\da \f = -\ii(\nda) \f$, since the ambient
connection is symmetric. Later on these notations and conventions for
the use of $\ii(\cdot)$ and $\e(\cdot)$ are also used for form
tractors, and related objects.

We write $\afl$ for the ambient {\em form Laplacian} $\da\ad+\ad\da$.
Using this with ambient form operators just introduced generates a
closed system of anti-commutators and commutators as given in Tables
\ref{anticomms} and \ref{comms}. (In fact the graded system is
isomorphically the Lie superalgebra $\fsl (2|1)$ and extends the $\fsl
(2)$ which played a role in \cite{GJMS}.  Some of the results below
could be rephrased as identities of $\fsl (2|1)$ representation
theory, but we have not taken that point of view. We also note that in 
\cite{HS}, which concerns powers
 of the ambient Dirac operator, the authors recover  
a 5-dimensional superalgebra isomorphic to the orthosymplectic
algebra ${\frak{osp}(2|1)}$.  This may be realised as a subalgebra
of ${\frak{sl}}(2|1)$.)  Note that the
relations in Table \ref{anticomms} are essentially just definitions and
standard identities. The relations in the table of commutators follow
from the anticommutator results,  $\ad \Q =2\X$,
(\ref{ndaX}), and the usual identities of exterior calculus on
pseudo-Riemannian manifolds.  In particular, they hold in all
dimensions and to all orders.
\begin{table}[ht]
\begin{center}
\begin{tabular}{c|c|c|c|c|}
\{\cdot,\cdot\} & $\ad$ &    $\da$ &    $\e(\X)$ &    $\ii(\X)$ \\
\hline
$\ad$   &   $0$       &    $\afl$   &    $0$         &   $\LX$ \\
\hline
$\da$ & $\afl$ & $0$ & $\LX^*$ & $0$ \\
\hline
$\e(\X)$ & $0$ & $\LX^*$ & $0$ & $\Q$  \\
\hline
$\ii(\X)$ & $\LX$ & $0$ & $\Q$ & $0$  \\
\hline
\end{tabular}\end{center}
\caption{Anticommutators $\{\mathfrak{g}_1,\mathfrak{g}_1\}$}
\label{anticomms}
\end{table}
\begin{table}[ht]
\begin{center}
\begin{tabular}{c|c|c|c|c|c|c|c|c|c|}
$[\cdot,\cdot]$ & $\ad$ & $\da$ & $\e(\X)$ & $\ii(\X)$ & & $\afl$ & $\LX$ & 
$\LX^*$ & $\Q$ \\
\hline
$\afl$ & $0$ & $0$ & $-2\ad$ & $2\da$ & & $0$ & $2\afl$ & $-2\afl$ & $-2\KX$ \\
\hline
$\LX$ & $0$ & $-2\da$ & $2\e(\X)$ & $0$ & & $-2\afl$ & $0$ & $0$ & $2\Q$ \\ 
\hline
$\LX^*$ & $2\ad$ & $0$ & $0$ & $-2\ii(\X)$ & & $2\afl$ & $0$ & $0$ & $-2\Q$ \\ 
\hline
$\Q$ & $-2\e(\X)$ & $2\ii(\X)$ & $0$ & $0$ & & $2\KX$ & $-2\Q$ & $2\Q$ & $0$ \\
\hline
\end{tabular}\end{center}
\caption{Commutators $[{\mathfrak{g}}_0,{\mathfrak{g}}_1]$ and
$[{\mathfrak{g}}_0,{\mathfrak{g}}_0]$, where $\KX:=\LX-\LX^*$}
\label{comms}
\end{table}

Now since each section $V\in \ct^k[w]$ is equivalent to
$\tilde{V}\in \act^k(w)|_\cq$, it follows that operators along $\cq$
that correspond directly to operators on $\ct^k[w]$ should not depend on
how $\tilde V$ is extended off $\cq$. We say a differential operator acts
{\em tangentially} along $\cq$, if $P\Q=\Q P'$ (or $[P,\Q]=\Q(P'-P)$)
for some operator $P'$, since then 
$$
P (\tilde V + \Q U) = P\tilde V + \Q P'U
$$
and so $P \tilde{V}|_\cq$ is independent of how $\tilde{V}$ is
extended off $\cq$.  Note that compositions of tangential operators
are tangential. If tangential operators are suitably homogeneous then
they descend to operators on on $ M$ and, since the ambient manifold
does not depend on any choice of metric from the conformal class, the
resulting operators are conformally invariant.  Of course they may
depend on choices involved in the ambient metric, in which case they
would fail to be natural. We will return to this point shortly.

Consider the form Laplacian. From the commutator table we have
$[\afl,\Q]=-2\KX$ where $\KX$ is a shorthand for $\LX-\LX^*$. Thus in
general $\afl$ is not tangential. However via a standard
pseudo-Riemannian identity and \nn{ndaX} one has 
$$
\KX=\LX-\LX^* = n+2\NX +2 ,
$$ and so $\KX$ acts as the zero operator on ambient forms homogeneous
of weight $-1-n/2$. If $U$ is homogeneous of weight $-1-n/2$ then $\Q
U$ is homogeneous of weight $1-n/2$, and so $\afl$ does act
tangentially on $\act^k(1-n/2)$. The form Laplacian is homogeneous of weight $-2$ 
in the sense that 
$[\LX,\afl]=-2\afl$ and so $\afl$ determines a well-defined operator 
$$
\afl: \act^k(1-n/2)|_\cq \to \act^k(-1-n/2)|_\cq ,
$$
which is clearly equivalent to an operator between tractor bundles, that we shall 
denote $\fb$:
$$
\fb: \ct^k[1-n/2]\to \ct^k[-1-n/2].
$$

This example generalises. On any ambient form field one has   
$$ 
[\afl^m,\Q]=\sum_{p=0}^{m-1}\afl^{m-1-p}[\afl,\Q]\afl^p =
-2\sum_{p=0}^{m-1}\afl^{m-1-p} \KX \afl^p. 
$$
From the homogeneity of $\afl$ it follows that, acting on 
$\act^k(w)$, the $p^{\underline{{\rm th}}}$ term on the right acts as
$-2[2(w-2p)+n+2]\afl^{m-1}$. Summing terms we get that 
$[\afl^m,\Q]$ acts
as $-2m(2w-2m+n+4)\afl^{m-1}$ on $\act^k(w)$ and so 
\begin{equation}\label{powerslap}
\afl^m:\act^k(m-n/2)\to\act^k(-m-n/2) \quad \mbox{is tangential}.
\end{equation}
 We write $\fl_m$ for the corresponding conformally invariant operator
on $M$: $\fl_m:\ct^k[m-n/2]\to\ct^k[-m-n/2].  $ Note that because of
the weights involved the formal adjoint operator maps between the same
spaces, $ \fl_m^*: \ct^k[m-n/2]\to\ct^k[-m-n/2]$.  It seems likely
that these operators agree. (For example they do if $k=0$
\cite{GrZ,FGrQ}.) However we do not need to investigate this since we
can simply work with the formally self-adjoint average of these
$$
 \tfrac{1}{2}(\fl_m+\fl_m^*)=:\fb_m:\ct^k[m-n/2]\to\ct^k[-m-n/2]
 \quad m\in 
\{0,1,2,\ldots\}.
$$
It is straightforward to verify that these have leading term $(-1)^m\Delta^m$.

The above conformally invariant powers of the Laplacian arise from
ambient operators which are tangential only for a specific
weight. There are also ambient operators which act tangentially on
forms without any assumptions of homogeneity.  From the tables we have
$[\KX,\Q]= [\LX-\LX^*,\Q]=4\Q$ and so $(\KX-4)\Q = \Q \KX$. Then
$[\d,\Q]=-2\ii(\X)$. On the other hand $[\ii(\X)\afl,\Q]=
\ii(\X)[\afl,\Q]= -2\ii(\X)\KX$. So 
$$
\ii(\afD):=-\da(\LX-\LX^*-4) + \ii(\X)\afl 
$$
satisfies
$$
[\ii(\afD),\Q]=-4\Q \d ,
$$
which shows that 
$
\ii(\afD)
$ acts tangentially on any ambient form. Similarly for $\e(\afD): =
\ad(\LX-\LX^* -4) +\e(\X)\afl$. These are each homogeneous of weight
$-1$ so for each $w\in {\Bbb R}$ we have  tangential operators
$$
\e(\afD):\act^k(w)\to\act^{k+1}(w-1),\ \
\ii(\afD):\act^k(w)\to\act^{k-1}(w-1).
$$ 
Via the identities of the tables and the others discussed, there
are many alternative ways to write these operators. For example we have 
\begin{equation}\label{altD}
\e(\afD)=
(\LX-\LX^*)\ad +\afl\e(\X)=(n+2\Euler+2) \ad+\afl \e(X). 
\end{equation}
The corresponding conformally invariant operators on form tractors
are denoted respectively $\e(\fD)$ and $\ii(\fD)$,
$$
\e(\fD):\ct^k[w]\to\ct^{k+1}[w-1],\ \
\ii(\fD):\ct^k[w]\to\ct^{k-1}[w-1].
$$ 
Using the tables it is straightforward to show that these satisfy
many surprising and useful useful identities, for example
 $$
\ii(\fD)\ii(\fD)=0,\ \ \e(\fD)\e(\fD) =0,\ \
\ii(\fD)\e(\fD)+\e(\fD)\ii(\fD)= 0 .
$$ 

The main key to our constructions in the next lecture is the
 result that, as an operators on $\ct^k[1+\ell-n/2]$, we have the remarkable identity
\begin{equation}\label{key2}
\fb_{\ell} \e(\fD) =\e(X) \fb_{\ell+1} ,
\end{equation}
which generalises significantly earlier known identities for low order
and the conformally flat case \cite{gosrni}.  Even at low orders,
verifying this using explicit formulae for the operators on $M$ would
be a daunting task. At the ambient
level this an almost trivial consequence of the fact that the form
Laplacian $\afl$ commutes with both $\ad$ and $\da$.  The relevant
result there is worthy of some emphasis so we write it as a
proposition.
\begin{proposition} \label{old-domino}
If $ V\in \act^{k}(\ell-n/2+1)$ and $ U\in \act^{k+1}(\ell-n/2)$ then
for $\ell =0,1,\cdots $ we have
$$
\afl^{\ell} \e(\afD)V =\e(\X) \afl^{\ell+1} V,  
\quad \ii(\afD)\afl^{\ell} U = \afl^{\ell+1}\ii(\X) U .
$$
Here $ \afl^0$ means 1.
\end{proposition}
\noindent{\bf Proof:} We will prove the first identity; the proof of the other
is similar.  
First observe that acting on any ambient form field, we have
$$
\afl^{\ell}(2\ell \ad+ \e(\X)\afl) 
=2\ell \afl^{\ell}\ad + \afl^{\ell}\e(\X)\afl 
= 2\ell \afl^{\ell}\ad + [\afl^{\ell},\e(\X)]\afl  + \e(\X) \afl^{\ell+1} .
$$
Now recall from the Tables \ref{anticomms} and \ref{comms}
that $ [\afl,\e(\X)]=-2\ad$, and that $\afl$ and $\ad$ commute.
Thus
$ [\afl^{\ell},\e(\X)]\afl=-2 \ell \afl^{\ell}\ad$, giving
$$
\afl^{\ell}(2\ell \ad+ \e(\X)\afl) = \e(\X) \afl^{\ell+1}.
$$
 On the other hand, from the definition of $ \e(\afD)$, we have that
$\e(\afD) V =(2\ell \ad +\e(\X)\afl) V$ for $ V\in
\act^{k-1}(\ell-n/2+1)$. $\qquad\Box$\\ 
Interpreting the proposition
down on the underlying conformal manifold $M$, the second display of
the proposition gives $\ii(\fD)\fl_{\ell} = \fl_{\ell+1}\ii(X) $ on
$\ct^{k+1}[\ell-n/2]$. The formal adjoint of this is $\fl^*_{\ell}
\e(\fD) =\e(X) \fl^*_{\ell+1}$ on $\ct^k[1+\ell-n/2]$, while the other
display of the proposition gives $\fl_{\ell} \e(\fD) =\e(X)
\fl_{\ell+1}$ on the same space. Thus combining these gives \nn{key2}.

\subsection{Naturality}

The construction of the ambient manifold does not depend on choosing a
particular metric from the conformal class. So it is conformally
invariant. However the manifold is not unique. For our purposes we
have fixed some choice of ambient metric and we must check that
operators finally obtained on $M$ via the ambient construction do not
depend on the choices made in arriving at our particular ambient
manifold. It clearly suffices to show that the operators are natural
for the underlying conformal structure on $M$ i.e., given by universal
formulae polynomial in the metric $g$ and its inverse, the Levi
connection, its curvature and covariant derivatives. In fact there is
an algorithm \cite{GoPet,deRham} for expressing the operators we
require here (and in fact a significantly larger class \cite{GoPet2}) in terms of
such a formula. Here we sketch some of the ingredients.

One basic idea behind this algorithm is to understand how certain
tractor fields and operators, the explicit formulae for which are
already known, turn up on the ambient manifold.  Closely related to
$\e(\afD)$ and $\ii(\afD)$ is the operator $\D:=\nda (n+2\NX-2) +\X
\al$. This acts tangentially on any ambient tensor field and since it
is homogeneous determines an operator between tractor bundles. It is
straightforward to verify that this resulting operator is exactly the
tractor-D operator which we met earlier in the proof of proposition
\ref{notstrong}.  Recall that for a metric from the conformal class,
is given by
$$
D_A V:=(n+2w-2)w Y_A V + (n+2w-2)Z_A{}^{a}\nabla_a V  -X_A( \nd_p\nd^p V+w\J V),
$$
where $ V$ is a section of any tractor bundle of weight $w$.
In these formulae $\nd$ means the coupled tractor--Levi-Civita connection.

Next  consider
the tractor curvature $K$, which is defined by
$[\nabla_a,\nd_b]V^A=K_{ab}{}^A{}_B V^B$ for $V\in \ct$. In terms
of a choice of metric from the conformal class, this has the formula
$$ 
K_{abCE}= Z_C{}^cZ_E{}^e
C_{abce}-4X_{[C}Z_{E]}{}^e\nd_{[a}\Rho_{b]e}.
$$
This may be inserted invariantly into the space $\ct^3\otimes
\ct^2$ by $K \mapsto {\bbX}^3\cdot K=\e(X)\bbZ^2\cdot K$, where for
taking the inner product implicit in the `$\cdot$' we view $K$ as a
2-form (taking values in $\ct^2$). Let us write $\Omega$ for
$\bbZ^2\cdot K$.  Then from the relationship of the tractor connection
to the ambient connection, as we outlined earlier, it follows easily
that the conformally invariant field $\e(X)\Omega$ is exactly the
tractor field corresponding to $\e(\X)\aR$.  Using this and the
ambient Bianchi identity one finds that, in dimensions other 4, the
tractor field on $M$ exactly corresponding to to the ambient curvature
$\aR$ is $W:=\frac{3}{(n-2)(n-4)}\ii(D)\e(X)\Omega$ which is readily
expanded to give an explicit formula for $\aR$ in terms of the
underlying conformal structure.

The difference between the ambient Bochner and form
 Laplacians is given by
$$
\afl=-\al-\aR\hash\hash ,
$$
(where $\al:=\nda^A\nda_A$.)
It follows immediately that, as operators on $k$-forms,
\begin{equation*}
\e(\afD) =\e(\D) -\e(\X)\aR\hash\hash,\ \
\ii(\afD) =\ii(\D) -\ii(\X)\aR\hash\hash ,
\end{equation*}
and so for the corresponding operators on $M$ we have 
$$\textstyle
\e(\fD) =\e(D)- \e(X) \Omega \hash \hash ,\ \
\ii(\fD) =\ii(D)- \ii(X) \Omega \hash\hash .
$$ 
So from the formulae for $D$ and $K$ we have naturality and explicit formulae for these operators.

For higher order operators the first step is to express the relevant
ambient operator as a composition of low order tangential operators
that we already understand in terms of natural operators on $M$.  Let
us consider, for example, the operator $ \fb_2: \ct^k[2-n/2]\to
\ct^k[-2-n/2].  $ This arises from $\afl^2$ on the ambient manifold.
 Since $\afl^2$ acts tangentially
on $\tilde{V}\in\act^k(2-n/2)$, we are free to chose an extension of $\tilde{V}$ off
$\cq$ that simplifies our calculations without affecting $\afl^2
\tilde{V}|_\cq$. For example, given $\tilde{V}|_\cq$, it is easily verified that we can
arrange that $\al \tilde{V}=O(\Q)$ (where $O(\Q)$ in meant in the sense of formal
power series). Then using that $\afl=-\al-\aR\hash\hash $ we have 
$$
\begin{array}{rl}
{\afl^2 \tilde{V} }
&=(\al+\aR\sharp\sharp)(\al+\aR\sharp\sharp)\tilde{V}
\vspace{2mm}
\\
&
=\al^2 \tilde{V}+(\al\aR)\sharp\sharp \tilde{V}+2(\nda_{|A|}\aR)\sharp\sharp\nda^{|A|}\tilde{V}+\aR\sharp\sharp(\aR\sharp\sharp \tilde{V})+O(\Q)
,
\end{array}
$$ where the bars around indices indicate that they are to be ignored
for the purposes of expanding the $\sharp$'s. From the Bianchi identity
and the Ricci flatness of the ambient metric (for simplicity we assume
that $n\notin \{4,6\}$) it follows that $(\al\aR)\sharp\sharp
\tilde{V}=-\frac{1}{2}\aR\sharp\sharp(\aR\sharp\sharp \tilde{V}) $ and so this
combines with the last term.  Similarly, since $\X_A\nda^A
V=\nda_{\miniX}\tilde{V}=(2-n/2)\tilde{V}$ and $\D_A \aR =(n-6)\nda_A
\aR -\X_A \al \aR$ we can replace $\nda_{|A|}\aR$ with
$\frac{1}{(n-6)}\D_A \aR$ provided we further adjust the coefficient of
the last term. On the other hand if we write $V\in \ct^k[2-n/2]$ for
the tractor field equivalent to $\tilde{V}|_\cq$ then it is clear that
$\aR\sharp\sharp(\aR\sharp\sharp \tilde{V})|_\cq$ is equivalent to $W\sharp\sharp(W\sharp\sharp V)$.
We have an explicit formula for $W$ and so we can deal with these
terms. Next note that since $\al \tilde{V}=O(\Q)$ it follows that
$2\nda^{A}\tilde{V}=2\nda^{A}\tilde{V}-\X^A\al \tilde{V} +O(\Q) =\D^A
\tilde{V} +O(\Q)$. Thus, modulo terms of the form $\aR\sharp\sharp(\aR\sharp\sharp
\tilde{V})$, and modulo $O(\Q)$, the term
$2(\nda_{|A|}\aR)\sharp\sharp\nda^{|A|}\tilde{V}$ is a multiple of
$(\D_{|A|}\aR)\sharp\sharp\D^{|A|}\tilde{V}$. The restriction of this to $\cq$
is equivalent to $(D_{|A|}W)\sharp\sharp D^{|A|} V$ which is a combination of
standard tractor objects for which we have explicit formulae. 

Finally there is the term $\al^2 \tilde{V}$.  From the relationship
between $\D$ and $D$ it follows easily that $\al$ acts tangentially on
ambient tensors homogeneous of weight $(1-n/2)$ and in this case
descends to the conformally invariant $\Box$ on tractor fields of
weight $(1-n/2)$ on $M$. The operator $\D$ lowers homogeneous weight
by 1 and so $\al \D$ acts tangentially on $\tilde{V}$ and $\al
\D\tilde{V}|_\cq$ is equivalent to $\Box D V$. Once again $\Box D$ is
a composition of operators for which we have explicit formulae. This
deals with $\al^2 \tilde{V}$ since $\al \D V= \al (2 \nda -\X\al)
\tilde{V} =-\X\al^2 \tilde{V} +2 [\al,\nda]\tilde{V}$ and expanding
the commutator here yields ambient curvature terms each of which is
either $O(\Q)$ due to the conditions on the ambient Ricci curvature or
can be dealt with by a minor variation of the ideas discussed above. 

More generally for the powers $\afl^m$ (as in \nn{powerslap}) that we
require it is straightforward to show inductively \cite{GoPet2} that
there is an algorithm for re-expressing each of these, modulo $O(\Q)$,
as an operator polynomial $\X$,$\D$, $\aR$, the ambient metric $\h$,
and its inverse $\h^{-1}$. A formula for the corresponding operator on
$M$ is then given by a tractor expression which is the same formula
except with $\X$,$\D$, $\aR$, $\h$, and $\h^{-1}$ formally replaced
by, respectively, $X$, $D$ $W$, $h$ and $h^{-1}$.

\section{Lecture 3 -- Proving the Theorem -- a sketch}

Henceforth we restrict ourselves to the setting where the underlying
conformal manifold is even dimensional.

Recall that the conformal gauge operator that we constructed for the
Maxwell operator arose from a single tractor operator, taking values in
a subbundle of the tractor bundle $\bT^{n/2-1}[-1-n/2]$ consisting of
elements of the form
$$
\left(\begin{array}{c}0\\
* \qquad 0 \\
*
\end{array}\right).
$$ Since the tractor-projectors $X$ and $Y$ are null and $\ii(X)Y=1$ it
follows at once from \nn{parts} that this subbundle may alternatively
be characterised as the elements of the form $\ii(X)F$. Let us write
$\bG_k$ (with section space denoted $\cg_k$) to denote the subbundle
of $\bT^k[k-n]$ consisting of elements of the form $\ii(X)F$ and write
$\bG^k$ (with section space $\cg^k$) for the quotient bundle which
pairs with this in conformal integrals. That is $\bG^k:=
\bT^k[-k]/\bV^k[-k]$ where $\bV^k[-k]=\Ker(\e(X):\bT^k[-k]\to
\bT^{k+1}[1-k])$ (equivalently $\bV^k[-k]$ is the subbundle of
$\bT^k[-k]$ consisting of elements of the form $\e(X)U$). Note that
from the composition series \nn{splitmax} for $\ct^k$ we see that 
the section spaces have composition series

\vspace{-23pt}

$$
\cbGuk\quad \mbox{and} \quad \cbGdk . 
$$ 

\vspace{9pt}

We will write $q_k:\ce^k\to \cg^k$ for the canonical conformally
invariant injection which, in a choice of conformal scale, is given by
$\bbZ^k\cdot$.  Similarly we write $q^k:\cg_k\to \ce_k$ for the
conformal map to the quotient which, in a conformal scale, is given by
$\bbZ\bullet$, where the bullet indicates the inner product on the
form tractor indices. Note that the composition $\ii(X) q_k$ vanishes,
and so also the dual operator $q^k\e(X)$.

In this language the $(\mbox{Maxwell operator},\mbox{ gauge})$ pair
appearing as the first term in the sum \nn{splitmax} is really a
single conformally invariant operator $\ce^{n/2-1}\to \cg_{n/2-1}$.
To generalise this we should obviously find an operator, $\ce^k\to
\cg_k$. However the symmetry of Maxwell detour complex suggests that
there might be a more general operator taking $\cg^k$ to $\cg_k$.
Given our observations above this is easily constructed. First observe
that $\e(X)$ gives a conformally invariant bundle injection
$\e(X):\bG^k\to \bT^{k+1}[1-k]$. The conformally invariant adjoint
operation is $ \ii(X): \bT^{k+1}[k-1-n]\to \bG_k $. Thus composing
$\fb_{\ell+1}$  fore and aft with these gives the
required operator,
$$
\big(\ii(X)\fb_{\ell+1}\e(X):=\bK_{\ell}\big): \cg^k \to \cg_k ,
$$ where $\ell:=n/2-k$.  Note that this is manifestly formally
self-adjoint since $\fb_{\ell+1}$ is formally self-adjoint. The
non-triviality is an easy consequence of the ellipticity of
$\fb_{\ell+1}$ and the classification of operators on the
conformal sphere (as discussed in lecture 2). Finally by the algorithm for constructing a formula
for these operators one finds they are natural for $k\in\{0,1,\cdots
,n/2\}$.

This at once gives candidates for the required long operators, viz.\ the full composition 
$L_k=q^k\bK_{\ell}q_k$ as in the diagram 
        \begin{center}\begin{picture}(436,60)(8,-30)

\put(99,0){$\ce^k$}        

\put(130,9){$q_k$}       
\put(123,4){\vector(1,0){20}}       

\put(155,24){$\ce^{k-1}$}
\put(153,12){\mbox{$\upl$}}
\put(155,0){$\ce^k$}        

\put(185,9){$\bK_\ell$}       
\put(179,4){\vector(1,0){20}}       

\put(212,0){$\ce_k$}      
\put(209,-13){$\upl$}
\put(212,-24){$\ce_{k-1}$}

\put(250,9){$q^k$}                     
\put(243,4){\vector(1,0){20}}          

\put(275,0){$\ce_k$.}        

\end{picture}\end{center}
\setlength{\unitlength}{1pt}
By construction the $L_k$ are formally self-adjoint and once again the
non-triviality follows from the non-triviality of the $\bK_\ell$ and
the classification of the operators between forms in the conformally
flat setting.  The candidates for extensions of the $L_k$ to
conformally invariant elliptic operators (as in part \IT{iv} of the
Theorem) are the operators obtained by simply omitting the final
projection, that is $\bL_{\ell}:=\bK_{\ell}q_k $.

Toward establishing that these operators have the desired properties
it is useful to make an observation related to the geometry of the
underlying constructions. This is that, since $2 \X=d\Q$, and
$\Q$ is a defining function for $\cq$ in $\aM$, it follows that
$\bG^k$ may be naturally identified with $\wedge^kT^* \cq/\sim$ where
$U_{p}\sim V_q$ if $U_{p}= \rho^s_* (V_q)$ for some $s\in {\Bbb R}_+$. This is an easy consequence of our
recovery of the tractor bundles from a similar quotient of $T\aM\cong T^*\aM$.
The
exterior derivative \underline{on $\cq$} preserves the subspace of
forms homogeneous with respect to  the canonical ${\Bbb
R}_+$ action $\rho(s)$ and so this determines an operator (for each
$k$) 
$$
\td: \cg^k \to \cg^{k+1} \quad \mbox{with formal adjoint} \quad \dt: \cg_{k+1}\to \cg_k.
$$
Under the described geometric interpretation of $\cg^k=\ce^{k-1}\lpl\ce^k$ the $\ce^k$-part
of the composition series arises from $\pi^*\ce^k$. Since exterior differentiation commutes 
with the pull-back it follows that 
$$ 
\td q_k =q_{k+1} d \quad \mbox{on $\ce^k$} \quad \mbox{and} \quad
q^k \dt = \d q^{k+1} \quad \mbox{on $\ce_{k+1}$},
$$
where the second result follows by taking the formal adjoint of the first.

More generally we get an operator $\td:\cg^k[w]\to \cg^{k+1}[w]$ given by
$$
\bbY^k
\cdot \alpha +\bbZ^k \cdot \m =
\left( \begin{array}{c}
\alpha\\
\m
\end{array}\right) \mapsto
\left( \begin{array}{c}
w\mu-\e(\nd)\alpha\\
\e(\nd)\m
\end{array}\right)
=
\bbY^{k+1}\cdot(w\mu-\e(\nd)\alpha)+\bbZ^{k+1}\cdot \e(\nd)\m ,
$$ and a formal adjoint $\dt$ for this.  This generalisation of $\td$
still arises from the exterior derivative on $\cq$, except now
restricted to appropriately homogeneous sections of
$\Lambda^kT^*\cq$. From these origins, or alternatively the explicit
formula displayed, it is clear that $\td^2=0$ (and hence also
$\dt^2=0$). Also, the operator $\td$ satisfies the anti-derivation rule 
$\td(\e(U)V)=\e(\td U) V+(-1)^k \e(U) \td V$ 
for $U$ in $\cg^k[w]$ and $V$ in any $\cg^{k'}[w']$. 

These operators turn up as factors in the components of $\bK_\ell$. Consider
$\e(X)\bK_{\ell}$. As a map on $\cg_k$ we have 
\begin{equation}\label{eXmeans}
\e(X): \bGdk \hspace*{4mm} \rightarrow \hspace*{4mm} \bGdkone \quad \mbox{by} \quad 
\left(\begin{array}{c}
u\\
v
\end{array}\right) \mapsto 
\left(\begin{array}{c}
0\\
u
\end{array}\right).
\end{equation}
 Now $\e(X)\bK_{\ell}= \e(X)\ii(X)\fb_{\ell+1}\e(X)$. Since $X$
is null, $\e(X)$ and $\ii(X)$ anticommute. Then by \nn{key2} we
have $ \e(X)\fb_{\ell+1}=\fb_\ell \e(\fD)$.  Thus
$\e(X)\bK_{\ell}=-\ii(X)\fb_\ell \e(\fD) \e(X)$. But recall $\e(\fD)
\e(X)$ arises from the ambient composition $\e(\afD) \e(\X) $ and by
\nn{altD} this is exactly $ (n+2\Euler +2)\ad \e(\X)= - (n+2\Euler
+2)\e(\X)\ad $, where finally we have used the anti-commutativity of
$\e(\X)$ and $\ad$.  Thus on $\cg^k$ (which is a quotient of
$\ct^k[-k]$) $\e(\fD) \e(X) = -(n-2k+2)\e(X)\td=-2(\ell+1)\e(X)\td
$. Hence overall we have $\e(X)\bK_{\ell}= 2(\ell+1)
\ii(X)\fb_{\ell}\e(X)\td= 2(\ell+1) \bK_{\ell-1} \td $.  By taking the 
formal adjoint there  
is a
corresponding  result for $\bK_{\ell+1} \ii(X) $, and we summarise these
surprising results in the following lemma which is central to the
subsequent discussions.
\begin{lemma} \label{key} 
As operators on $\cg^k$ we have
$$
\e(X)\bK_{\ell} = 2(\ell+1) \bK_{\ell-1} \td\quad 
\mbox{and}\quad \bK_{\ell+1} \ii(X) = 2(\ell+2)\dt  \bK_\ell  .
$$
\end{lemma}

We are now ready to construct the operators of the Theorem. We start
by looking at $\bL_\ell=\bK_\ell q_k$. Note that $\ii(X)$ and $\e(Y)$
are well defined on $\cg^k$. Recall that on $\ct^k$, and therefore
also on $\cg^k$,
$$
\ii(X)\e(Y)+\e(Y)\ii(X)=h(X,Y)=1.
$$
Using this, and since $\ii(X)q_k=0$, we have $\bK_\ell q_k=
\bK_\ell \ii(X)\e(Y)q_k$. Thus from the lemma we have 
\begin{equation}\label{bLexp}
\bK_\ell q_k = 2(\ell+1)\dt \bK_{\ell-1} \e(Y)q_k.    
\end{equation}
The claim is that this single conformally invariant operator gives the
conformally invariant elliptic system $( \d Q_{k+1} d,~ \d Q_k)$ of part
\IT{iv} of the theorem. 

$\bK_\ell q_k $ takes values in $\cg_k=\ce_k\lpl \ce_{k-1}$. As
observed above, to obtain the component with range $\ce_k$ we project
with $q^k$.  Applying this to both sides of the last display we have
$L_k:= q^k \bK_\ell q_k = 2(\ell+1) q^k \dt \bK_{\ell-1} \e(Y)q_k $.
Then recall $q^k \dt = \d q^{k+1}$ so 
$$ L_k= 2(\ell+1) \d q^{k+1}\bK_{\ell-1} \e(Y)q_k $$ which establishes
that $L_k$ has $\d$ as a left factor. Continuing on to show that it also has $d$ as
a right factor involves similar arguments. First observe that $\ii(Y) \e(X)+\e(X)\ii(Y)$
is well-defined and acts as the identity on $\cg_{k+1}$ and recall
$q^{k+1}\e(X)=0$. So we may insert $\ii(Y) \e(X)$ to obtain 
$$ 
L_k= 2(\ell+1) \d q^{k+1}\ii(Y) \e(X) \bK_{\ell-1} \e(Y)q_k= 4 \ell
(\ell+1) \d q^{k+1}\ii(Y) \bK_{\ell-2} \td \e(Y)q_k ,
$$
where to obtain the last expression we have used the Lemma to
exchange $\e(X) \bK_{\ell-1}$ with a multiple of $ \bK_{\ell-2} \td$.
On the other hand it is easily verified that, in a conformal scale
$\si\in \ce[1]$, the corresponding $\e(Y)$ agrees with $\e(\si^{-1}\td
\si)$. From the anti-derivation rule for $\td$, $\td^2=0$, and that
the log-derivative $\si^{-1}\td \si $ is annihilated by $\td$, it
follows that $\td$ anti-commutes with $\e(Y)$. Finally we have already
noted that $\td q_k =q_{k+1} d$ and so by this general construction we
obtain, for each $k$, long operators which factor through $\d$ and
$d$,
$$
L_k= - 4 \ell
(\ell+1) \d q^{k+1}\ii(Y) \bK_{\ell-2} \e(Y) q_{k+1} d . 
$$ 
This completes part \IT{iii} since from their non-triviality it
follows easily that each $L_k$ has leading term $(\d d)^{n/2-k}$ (at
least up to a constant multiple).

This result for the $L_k$ suggests 
\begin{equation}\label{M4Q}
M_{k}:=q^{k}\ii(Y) \bK_{\ell-1} \e(Y) q_{k}
\end{equation} 
(or some
multiple thereof) is a candidate for $Q_{k}$.

Next we examine the $\ce_{k-1}$-component of \nn{bLexp}.  In terms of
the projectors from lecture 1 this is the coefficient of
$\bbX^{k}=\e(X)\bbZ^{k-1}$.  Since $q^{k-1}$ is given, in a choice of
conformal scale, by $\bbZ^{k-1}\bullet$ it follows that
$q^{k-1}\ii(Y)\bK_\ell q_k $ is the $\ce_{k-1}$-component of $\bK_\ell
q_k$.  Composing $q^{k-1}\ii(Y) $ with the right-hand-side of
\nn{bLexp} brings us to $2(\ell+1)q^{k-1}\ii(Y) \dt \bK_{\ell-1} \e(Y)q_k
$. That $\{ \ii(Y), \dt\}$ vanishes on $\cg_{k+1}$ is just the formal
adjoint of the result for $\{\td, \e(Y) \}$ on $\cg^{k-1}$. We have
already that $q^{k-1} \dt =\d q^k$ and so
$$q^{k-1}\ii(Y)\bK_\ell q_k= - 2(\ell+1)\d q^{k}\ii(Y)\bK_{\ell-1}\e(Y) q_k =-2(\ell+1)\d M_{k} .$$

Summarising then, we have that in a choice of conformal scale 
$$
\bK_{\ell}q_k: \ce^k\to \bGdk\quad  \mbox{is given by} \quad u\mapsto 
\left(\begin{array}{c}
\d M_{k+1} du\\
\d M_k u
\end{array}\right), \hspace*{2mm} k=1,\cdots n/2-1,
$$ 
where $M_k$ is defined by \nn{M4Q} above and we have ignored the
details of non-zero constants.  Note that by construction $
\bK_{\ell}q_k$ is conformally invariant.  Therefore $\d M_k $ is
conformally invariant on the null space of $\d M_{k+1} d$. This is
more fundamental than the transformation formula. Nevertheless, 
observe that the transformation formula claimed in part \IT{iv} is now
immediate from the invariance of $ \bK_{\ell}q_k$ and the conformal
transformation formula $\widehat \bbZ=\bbZ+\e(\Up)\bbX$.
Thus part  \IT{iv} is proved.

For Part \IT{ii} of the Theorem we should examine the conformal
transformation law of $M_k =q^k \ii(Y) \bK_{\ell-1} \e(Y) q_{k}$. This
looks potentially complicated as neither $\ii(Y)$ nor $\e(Y)$ is
conformally invariant. Let us simplify initially and consider just
part of this viz.\ $\bK_{\ell-1} \e(Y) q_{k}$. Although $\bK_{\ell-1}$
maps $\cg^{k+1}\to \cg_{k+1}$ we will show that the image of
$\bK_{\ell-1} \e(Y) q_{k}:\cC^k\to \cg_{k+1}$ lies entirely in the
$\ce_{k}$ part of
$$
\cbGdkone .
$$

\vspace{10pt}

Suppose then that $\phi$ is a closed $k$-form and consider 
$\e(X)\bK_{\ell-1} \e(Y) q_{k} \phi $. By the lemma, and the commutation rules observed above, 
this gives (a constant multiple of)
$$ \bK_{\ell-2} \td \e(Y) q_k \phi=-\bK_{\ell-2} \e(Y) \td q_k \phi=
-\bK_{\ell-2} \e(Y) q_k d \phi= 0.
$$
From \nn{eXmeans} this exactly proves that in a choice of scale we have 
$\bK_{\ell-1} \e(Y) q_{k} \phi = \bbX^{k+1}\cdot(\tilde{M}_k) $, or in the matrix notation, 
$$
[\bK_{\ell-1} \e(Y) q_{k} \phi ]_g =
\left(\begin{array}{c}
0\\
\tilde{M}_{k} \phi
\end{array}\right)
$$ for some operator $\tilde{M}_k:\cC^k \to \ce_k$.  Recall one
recovers the coefficient of $\bbX^{k+1}$ by left composing with
$q^k\ii(Y)$ and so $\tilde{M}_k\phi= q^k\ii(Y)\bK_{\ell-1} \e(Y) q_{k} \phi
=M_k \phi $.  How does this formula transform conformally?  For the
purposes of this calculation we may take  $Y$ to be the metric dependent
section in $\cg^1[1]$ given by $Y=\si^{-1}\td \si$ 
 where $\si\in \ce[1]$ is the conformal scale,
that is $g=\si^{-2}\bg$. Then conformal rescaling
$$
g\mapsto \hat{g}= e^{2\omega}g, \quad \mbox{ corresponds to } 
\quad \si\mapsto\hat{\si}= e^{-\omega}\si,
$$
whence 
$$
Y\mapsto \widehat{Y}=Y-\td \omega .
$$
From this we have 
$$
\bK_{\ell-1} \e(Y) q_{k} \phi \mapsto \bK_{\ell-1} \e(\widehat{Y}) q_{k} \phi
= \bK_{\ell-1} \e(Y) q_{k} \phi - \bK_{\ell-1} \e(\Up) q_{k} \phi
$$ where, as usual $\Up:=d\omega$.  
 Now observe that, by the anti-derivation rule, $\td q_k \omega
\phi = \td \omega q_k\phi = \e (\Up) q_k \phi$, since $\td q_k=
q_{k+1} d$ and $d\phi=0$. Thus the conformal variational term in the
display may be written as the composition $\bK_{\ell-1} \td q_k \omega
\phi $. Using the Lemma we have $2(\ell+1)\bK_{\ell-1} \td =
\e(X)\bK_\ell$ and so the conformal variation term is $\e(X)\bK_\ell
q_k \omega \phi $ -- at least if we ignore the division by
$2(\ell+1)$.  Clearly this is also annihilated by the left action of
$\e(X)$ and so is of the form $\bbX^{k+1}\cdot \tilde{L}_k \omega \phi
$,
or 
$$
 \left(\begin{array}{c}
0\\
\tilde{L}_{k} \omega \phi
\end{array}\right) 
$$
for some operator $\tilde{L}_k:\ce^k\to\ce_k $. Once again we
recover the coefficient by composing on the left with $q^k\ii(Y)$ to
obtain $q^k\ii(Y)\e(X)\bK_\ell q_k \omega \phi = q^k\bK_\ell q_k \omega
\phi$ (since $\{ \ii(Y),\e(X)\}=1$ and $q^k\e(X)=0$) which we recognise
as $L_k \omega \phi$ (i.e.\ $\tilde{L}_k=L_k$).
So in summary we have shown the conformal transformation
$$
[\bK_{\ell-1} \e(Y) q_{k} \phi ]_g =
\left(\begin{array}{c}
0\\
M_{k} \phi
\end{array}\right) \mapsto [\bK_{\ell-1} \e(\widehat{Y}) q_{k} \phi ]_g = 
\left(\begin{array}{c}
0\\
M_{k} -\frac{1}{2(\ell+1)}L_k \omega \phi
\end{array}\right).
$$ 
Thus we exactly recover the result claimed in part \IT{ii} if we take 
$$
Q_k=(n+2)n\cdots (n-2k+2)M_k.
$$

\vspace{3mm}

\noindent {\bf Recovering Branson's Q-curvature}. So far in our
construction of the operators $\bK_\ell$, we have only used
$\afl^\ell$ on $\act^{n/2-\ell}(\ell-n/2)$ whereas as we have observed
already this operator is tangential on $\act^k(\ell-n/2)$ for any
$k$. The upshot of the latter observation is that for each $\ell$,
$\bK_\ell$ generalises to give an operator
$$
\bK_\ell:\cg^k[w]\to \cg_k[-w] ,
$$
where $w=\ell-n/2+k$. Thus we get order $2\ell$ conformally invariant
 operators between weighted forms
$$
(L_\ell^k:=q^k\bK_\ell q_k): \ce^k[w]\to \ce_k[w]. 
$$
These are natural for $\ell\leq n/2-1$, and also for $\ell=n/2$ if $k=0$. 
In this generalised setting Lemma \ref{key} still holds in the
sense that for example on $\cg^k[w]$ we have
$$
\e(X)\bK_\ell=2(\ell+1)\bK_{\ell-1}\tilde{d}.
$$

Thus for $L_\ell^k=q^k\bK_\ell q_k =q^k\ii(Y)\e(X)\bK_\ell q_k$ we obtain
the alternative formula
$$
L_\ell^k= 2(\ell+1)q^k \ii(Y)\bK_{\ell-1}\tilde{d}q_k.
$$ 
Now for $\mu\in \ce^k$, and $\si\in \ce[1]$ a conformal scale, we
have $\si^w \mu\in\ce^k[w]$.  We will apply $L^k_\ell$ to this. 
First note that in terms of the splittings for the conformal scale we
have
$$
q_k\si^w \mu =  \bbZ^k \cdot (\si^w \m). 
$$
So, from the explicit formula for $\td$, we have  
$$
\begin{array}{rcl}
\td q_k \si^w \mu &=&  w \si^w \bbY^{k+1}\cdot \mu+ \si^w \bbZ^{k+1}\cdot \e(\nd)\m\\ 
&=& w\si^w \e(Y)\bbZ^{k}\cdot \mu + \si^w \bbZ^{k+1}\cdot d \m ,\\
&=& w\si^w \e(Y)q_{k}\mu + \si^w q_{k+1} d \m ,
\end{array}
$$
where we used that the Levi Civita connection for $\si$ annihilates
$\si$. Using this again we obtain
$$
L_\ell^k \si^w \mu = 2 w (\ell+1) \si^w q^k \ii(Y)\bK_{\ell-1} \e(Y)q_{k} \mu+
2(\ell+1)\si^wq^k \ii(Y)\bK_{\ell-1} q_{k+1} d \m .
$$ 
If $\mu$ is closed then the second term vanishes and taking the
coefficient of $w$ and setting, in this, $\ell=n/2$ (i.e., $w=0$)
yields, up to a multiple, the $Q_k$ operator on $\mu$. In particular
if $k=0$ and we take $\mu=1$ we obtain $Q_0 1$. But in this case this
construction is exactly recovering the Q-curvature according to
Branson's original definition by dimensional continuation, since by
construction the operators $L^0_\ell$ agree with the GJMS operators
\cite{GJMS} on densities. 

We have not only shown that $Q_0 1$ is the usual Q-curvature but also
that the $Q_k$ operators, defined earlier without the use of
dimensional continuation, also arise from the analogous dimensional
continuation argument but now applied to the operators $L^k_\ell$
between form densities. Implicitly we are using that the operators
$\bK_\ell$ are given by universal formulae, polynomial in the
dimension, and with coefficients in terms of a stable basis of
Riemannian invariants. This is an immediate consequence of the tractor
formulae for the operators $\fb_\ell$ since the basic tractor objects
(such as the tractor-D operator and the $W$ tractor are given by
formulae in this form).

It is easily shown too \cite{deRham} that defining the Q-curvature to be $ Q_0 1$ 
with $Q_0$ defined by \nn{M4Q} is equivalent to  the definitions of both 
\cite{FeffH} and \cite{GoPet}.

\end{document}